\documentclass[pdflatex,sn-mathphys-num]{sn-jnl}
\usepackage{etoolbox}
\usepackage{paralist}
\usepackage{graphicx}
\usepackage{mathrsfs}
\usepackage{caption}
\usepackage{subcaption}
\usepackage{amssymb}
\usepackage{amsmath}
\usepackage{booktabs}
\usepackage{xcolor} 
\usepackage{multirow}
\usepackage{makecell}
\graphicspath{ {./} }
\usepackage{scrextend}
\usepackage{fancyhdr}
\usepackage{graphicx}
\newtheorem{lemma}{Lemma}
\newtheorem{theorem}{Theorem}
\usepackage{xcolor,colortbl}
\definecolor{Red}{RGB}{230,4,0}
\definecolor{Green}{RGB}{0,175,80}
\definecolor{Yellow}{RGB}{255,192,0}
\usepackage{subcaption}
\usepackage{algorithm2e}
\usepackage{lineno}
\usepackage[normalem]{ulem}
\usepackage{setspace}
\theoremstyle{thmstyletwo}
\newtheorem{remark}{Remark}[section]

\SetKwInput{KwInput}{Input}
\SetKwInput{KwOutput}{Output}

\def\equationautorefname#1#2\null{%
  Eq.#1(#2\null)%
}

\theoremstyle{thmstyleone}
\newtheorem{corollary}{Corollary}

\theoremstyle{thmstylethree}
\newtheorem{definition}{Definition}

\begin{document}

\title[Article Title]{Blend-to-zero operators for smooth transition functions}

\author[1]{\fnm{Iv\'{a}n} \sur{M\'{e}ndez-Cruz}}\email{ivan.mendez\_cruz@ens-paris-saclay.fr}

\author*[1]{\fnm{Faisal} \sur{Amlani}}\email{faisal.amlani@ens-paris-saclay.fr}

\affil[1]{ \orgname{Université Paris-Saclay, CentraleSupélec, ENS Paris-Saclay, CNRS, LMPS - Laboratoire de Mécanique Paris-Saclay}, \city{Gif-sur-Yvette}, \country{France}}








\abstract{Motivated by existing blend-to-zero techniques, a formal framework is developed for defining and constructing blend-to-zero operators on closed intervals for the generation of sufficiently smooth transitions between functions. Such transitions are first formulated as a two-point Hermite-type interpolation that is not necessarily polynomial. It is shown that, in the polynomial case, the corresponding interpolant can be explicitly represented in terms of the regularized incomplete Beta-function. This representation is then used to generate linear blend-to-zero operators. Following this, additional blend-to-zero operators are constructed by considering the algebraic and geometric properties of functions with sufficiently flat ends (e.g., smooth staircase functions and smooth step functions). Finally, explicit formulas for a family of trigonometric smooth step functions are provided, and these functions are shown to be related to certain higher-order two-point boundary value problems.}




\maketitle
\section{Introduction}\label{sec:intro}

Smooth transitions between functions are essential for generating partitions of unity \cite{Lee2012,Tu2011}, with applications ranging from several processes of computer graphics and computer-aided design \cite{Ebert2003,Hartmann2001,Li2004} to deep learning \cite{Shamir2020}. Such transitions can be generated as convex combinations of functions by using blending functions as weights \cite{Laksa2022}. Inspired by blending techniques, smooth unit step functions, and smooth staircase constructions \cite{Hartmann2001,Li2004,Laksa2022,Leopold2020}, as well as by related blend-to-zero constructions in certain periodic extension methodologies that are employed for high-order PDE solvers \cite{amlani2016,fontana2022,melkior2025,nainwal2025modified}, this contribution introduces and formalizes {\em blend-to-zero operators} on ${C}^{n}[a,b]$ for the generation of sufficiently smooth transitions between functions. Such operators flatten functions at one end while preserving values and higher-order derivatives at the other. This work presents explicit constructions and formulations of these operators and their relationship to smooth transition functions. 

The paper is organized as follows. \autoref{sec:transition} formally defines smooth transition functions and blend-to-zero operators, and formulates the transition problem as a two-point Hermite-type interpolation problem at the endpoints (not necessarily polynomial). It is also shown that smooth transitions can be generated by adding leftward and rightward blend-to-zero components. \autoref{sec:hermite} then considers the associated polynomial two-point Hermite interpolation problem, where we show that the corresponding interpolant can be represented explicitly in terms of the regularized incomplete Beta-function (this representation yields linear blend-to-zero operators). \autoref{sec:staircases} studies functions with flat ends, including smooth staircase functions and smooth step functions.  In particular, we show that: functions with flat ends form a linear space closed under multiplication; smooth step functions are closed under both multiplication and composition; and symmetry about the midpoint provides sufficient conditions to generate smooth step functions. \autoref{sec:blendop} constructs further linear blend-to-zero operators from functions with flat ends and, in particular, from smooth step functions. Finally, \autoref{sec:trigstep} provides an alternative derivation of the trigonometric blending functions of \cite{Olofsen2019}, presenting explicit formulas for the resulting trigonometric step functions and showing that they are related to certain higher-order two-point boundary value problems.

\section{Smooth transition functions and blend-to-zero operators}
\label{sec:transition}

Denote by $\mathbb N_0=\{0,1,2,\ldots\}$ the set of non-negative integers, and let
$\overline{\mathbb N}_0=\mathbb N_0\cup\{\infty\}$. For finite
$j\in\mathbb N_0$, denote $D^j=d^j/dx^j$, and denote by $I$ the identity
operator on $C^n[a,b]$. Throughout this work, when an order is equal to $\infty$, the corresponding
conditions are understood to hold for every finite derivative order. For compactness of notation, functions may also appear without their arguments (i.e., $f=f(x)$).

\begin{definition}
  Let $a,a_{0},b_{0},b\in\mathbb{R}$ be such that $a<a_{0}<b_{0}<b$, and 
  let $\ell,m,n,r\in\overline{\mathbb N}_0$ be such that $\ell\leq m$ and $r\leq n$.
  Let $f\in{C}^{m}[a,b_{0}]$ and let $g\in{C}^{n}[a_{0},b]$. 
  A {\em smooth transition function of orders $\ell$ and $r$ from $f$ to $g$ over 
  $[a_{0},b_{0}]$} is a function $h\in{C}^{p}[a,b]$ such that:
   \begin{enumerate}[(i)]
    \item $h(x)=f(x)$ for $a\leq x\leq a_{0}$;
    \item $h(x)=g(x)$ for $b_{0}\leq x\leq b$;
    \item $h$ is $\ell$-times continuously differentiable at $x=a_{0}$;
    \item $h\in{C}^{q}(a_{0},b_{0})$;
    \item $h$ is $r$-times continuously differentiable at $x=b_{0}$, 
   \end{enumerate}
  where $p=\min\{\ell,r\}$ and $q=\max\{\ell,r\}$.
  A smooth transition function has order $r$ if its two orders are equal to $r$,
  and it has arbitrary order if its two orders are equal to $\infty$. 
\end{definition}

\begin{figure}[htbp]
	\centering
	\includegraphics[width=.5\textwidth]{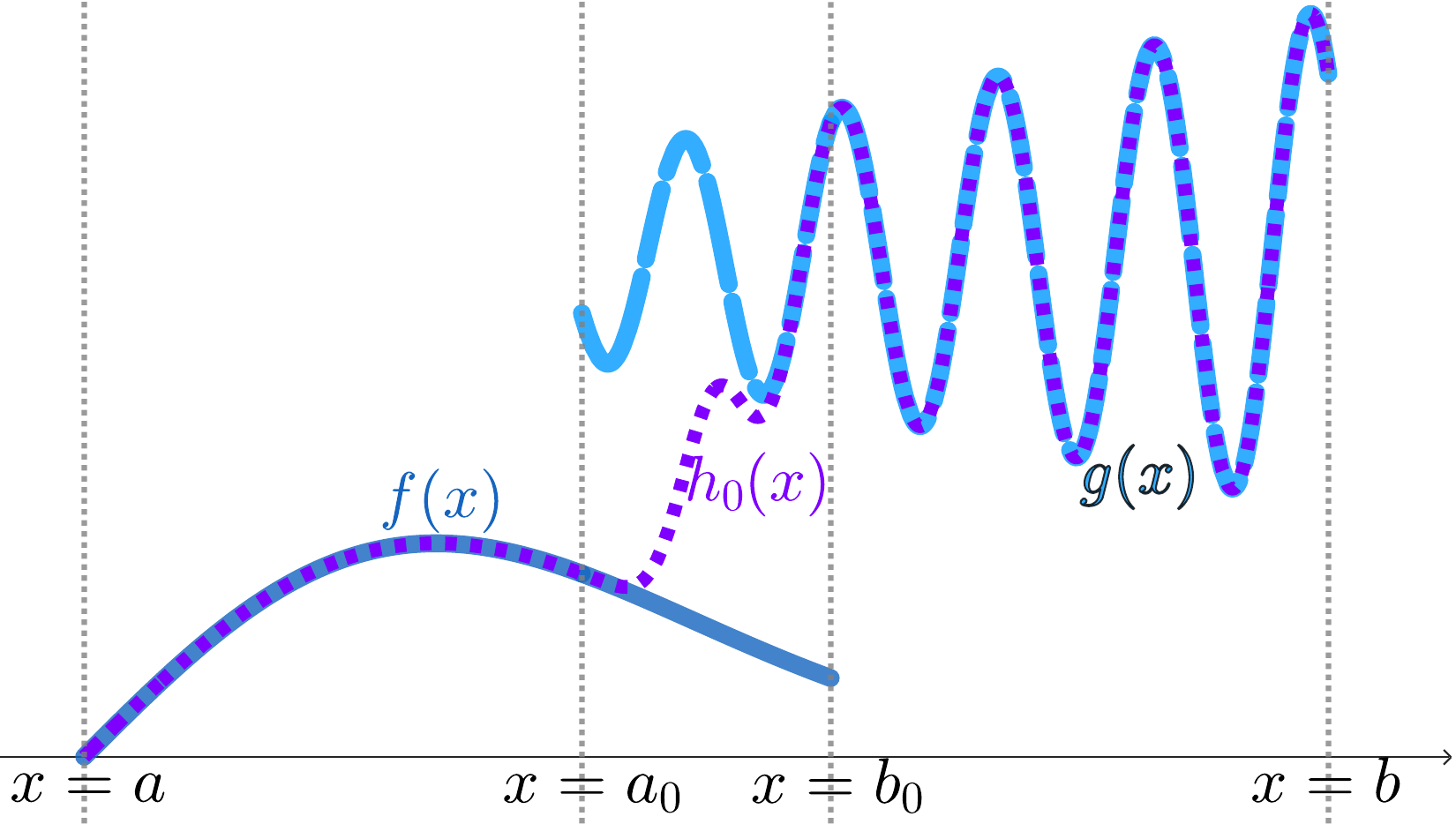}
  \caption{Illustration of a smooth transition function $h:[a,b]\to\mathbb R$ from $f:[a,b_0]\to\mathbb R$ to $g:[a_0,b]\to\mathbb R$. The transition function is shown as the purple dashed line, and its restriction to $[a_0,b_0]$ is denoted by $h_0$.}
	\label{fig:transition}
\end{figure}

\begin{remark}\label{rmk:interpol}
 A smooth transition function $h$ of orders $\ell$ 
  and $r$ from $f$ to $g$ over $[a_{0},b_{0}]$ satisfies the following:
\begin{alignat}{3}
  h_0^{(j)}(a_{0}) & = & f^{(j)}(a_{0}), && \quad j=0,\dots,\ell, \label{eq:fa}\\
  h_0^{(k)}(b_{0}) & = & g^{(k)}(b_{0}), && \quad k=0,\dots,r, \label{eq:gb}
\end{alignat}
where   $h_0$ is the restriction of $h$ to the interval $[a_0,b_0]$, i.e., $h_0 = h|_{[a_0,b_0]}$.
Conversely, if $h_{0}\in C^{q}[a_{0},b_{0}]$ satisfies the endpoint conditions in \autoref{eq:fa} and \autoref{eq:gb}, then the
piecewise function given by
  \begin{equation}\label{eq:pieceh}
    h(x) := 
    \begin{cases}
      f(x), & a\leq x < a_{0}; \\
      h_{0}(x), & a_{0}\leq x\leq b_{0}; \\
      g(x), & b_{0} < x \leq b \\
    \end{cases}
  \end{equation}
  is a smooth transition function of orders $\ell$ and $r$ from $f$ to $g$ over 
  $[a_{0},b_{0}]$. This construction is illustrated in \autoref{fig:transition}.
\end{remark}

\begin{remark}\label{rmk:sumtrans}
  Let $f_{0},g_{0}\in{C}^{q}[a_{0},b_{0}]$ be such that: 
  \begin{alignat}{8}
    f_0^{(j)}(a_{0}) & = & \ f^{(j)}(a_{0}), && \quad j=0,\dots,\ell,
    && \quad f_0^{(k)}(b_{0}) & = & 0, \qquad\ \ && \quad k=0,\dots,r, \label{eq:f0} \\
    g_0^{(j)}(a_{0}) & = & 0, \qquad\ \ && \quad j=0,\dots,\ell, 
    && \quad g_0^{(k)}(b_{0}) & = & \ g^{(k)}(b_{0}), && \quad k=0,\dots,r. \label{eq:g0}
  \end{alignat}
  By \autoref{rmk:interpol}, we have that $f_0$ generates a smooth transition function 
  of orders $\ell$ and $r$ from $f$ to the zero function over $[a_{0},b_{0}]$, while $g_{0}$ 
  generates a smooth transition function of the same orders from the zero function to $g$ over 
  $[a_{0},b_{0}]$. Furthermore, $h_0 = f_{0}+g_{0}$ satisfies \autoref{eq:fa} and \autoref{eq:gb}, and consequently it generates a smooth transition function of orders $\ell$ and $r$ from $f$ to $g$ over $[a_{0},b_{0}]$.
\end{remark}

In view of \autoref{rmk:sumtrans}, we now construct operators
that flatten functions at one endpoint while preserving their values and 
higher-order derivatives at the opposite endpoint.

\begin{definition}
  Let $\ell,q,r\in\overline{\mathbb N}_0$ be such that $\ell,r\le q$. An operator 
  $\mathcal{B}_{R}:{C}^{q}[a_{0},b_{0}]\to{C}^{q}[a_{0},b_{0}]$ 
  is called a {\em rightward blend-to-zero operator of orders $\ell$ and $r$} 
  if, for all $f\in{C}^{q}[a_{0},b_{0}]$, the function $f_{0}=\mathcal{B}_{R}(f)$ 
  satisfies \autoref{eq:f0}. An operator 
  $\mathcal{B}_{L}:{C}^{q}[a_{0},b_{0}]\to{C}^{q}[a_{0},b_{0}]$ 
  is called a {\em leftward blend-to-zero operator of orders $\ell$ and $r$} 
  if, for all $g\in{C}^{q}[a_{0},b_{0}]$, the function $g_{0}=\mathcal{B}_{L}(g)$ 
  satisfies \autoref{eq:g0}. A leftward or rightward blend-to-zero operator has order $r\in\mathbb N_0$ if its two orders are equal to $r$. It has arbitrary order if its two orders are equal to $\infty$.
\end{definition}

\begin{figure}[htbp]
	\centering
	\includegraphics[width=.4\textwidth]{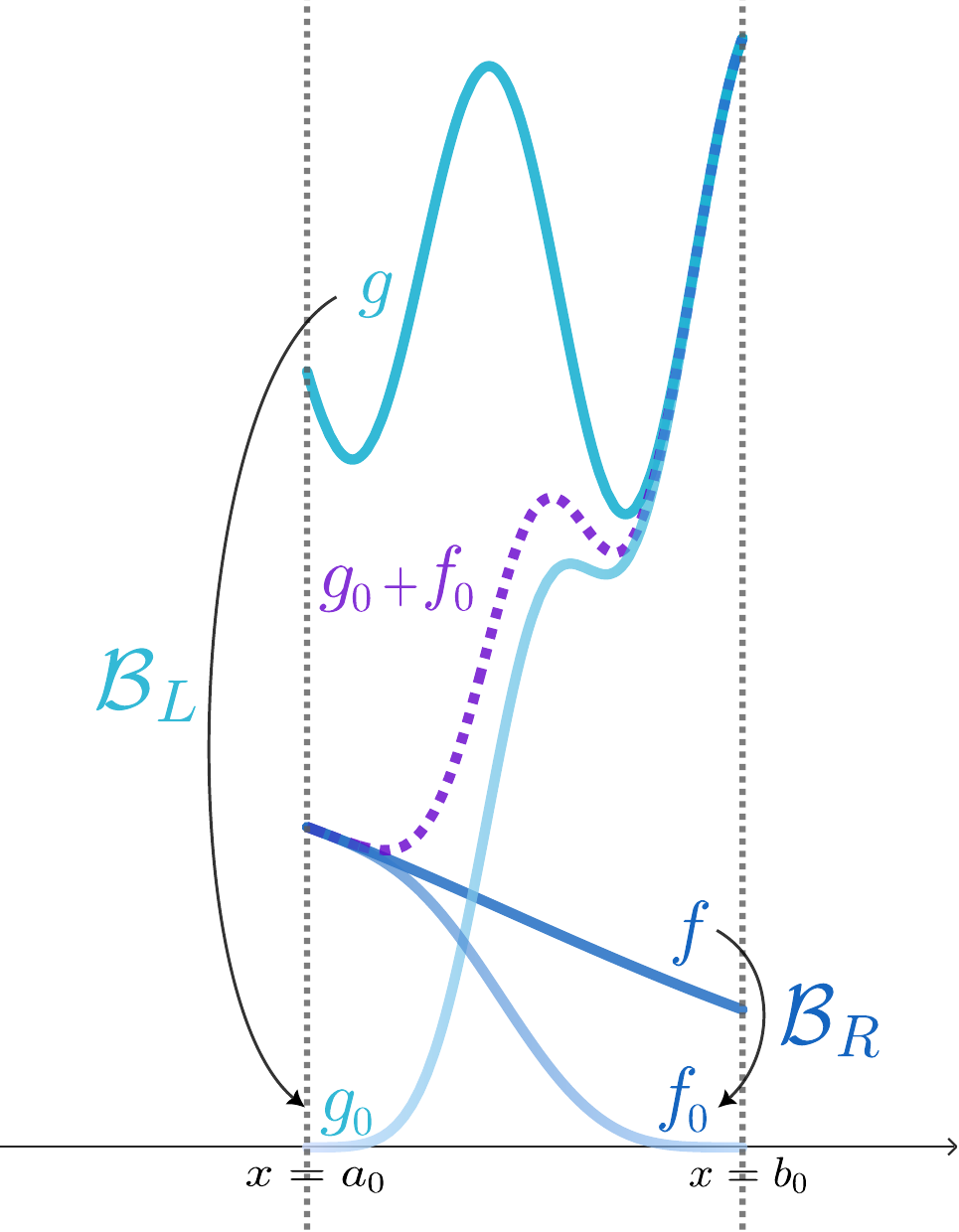}
	\caption{Illustration of a smooth transition from $f$ to $g$ over $[a_{0},b_{0}]$ 
obtained as the sum of $f_{0}=\mathcal{B}_{R}\!\left(f|_{[a_0,b_0]}\right)$ and  
  $g_{0}=\mathcal{B}_{L}\!\left(g|_{[a_0,b_0]}\right)$. The leftward blend-to-zero operator 
$\mathcal{B}_{L}$ flattens $g$ at $x=a_{0}$ and preserves the endpoint value and 
derivatives of $g$ at $x=b_{0}$, while the rightward blend-to-zero operator 
$\mathcal{B}_{R}$ preserves the endpoint value and derivatives of $f$ at $x=a_{0}$ 
and flattens $f$ at $x=b_{0}$.}
	\label{fig:sumblends}
\end{figure}

Smooth transition functions can be generated by adding blend-to-zero operators as illustrated in \autoref{fig:sumblends}. The following provides sufficient conditions for generating such functions.

\begin{theorem}\label{thm:sumblendop}
  Let $\mathcal{B}_{L}$ be a leftward blend-to-zero operator
  on $C^{q}[a_{0},b_{0}]$ of orders $\ell$ and $r$, and let 
  $\mathcal{B}_{R}$ be a rightward blend-to-zero operator
  on $C^{q}[a_{0},b_{0}]$ of the same orders. 
  Let $f\in C^{q}[a,b_{0}]$ and $g\in C^{q}[a_{0},b]$.  Define $f_{0}=\mathcal{B}_{R}\!\left(f|_{[a_0,b_0]}\right)$, 
  $g_{0}=\mathcal{B}_{L}\!\left(g|_{[a_0,b_0]}\right)$,
  and let $h_{0}=f_{0}+g_{0}$. 
  Then $h_0$ belongs to $C^{q}[a_{0},b_{0}]$ and satisfies \autoref{eq:fa} and \autoref{eq:gb}.
  Furthermore, the piecewise function $h$ defined by \autoref{eq:pieceh} is a smooth transition function of orders $\ell$ and $r$ from $f$ to $g$ over $[a_0,b_0]$.
\end{theorem}

\begin{proof}
 Since both $\mathcal{B}_{L}$ and $\mathcal{B}_{R}$
  are operators from $C^{q}[a_{0},b_{0}]$ into itself,
  we have that $f_{0},g_{0}\in C^{q}[a_{0},b_{0}]$.
  Hence $h_0\in C^{q}[a_{0},b_{0}]$.  Since $\mathcal{B}_{R}$ is a rightward blend-to-zero operator of orders 
  $\ell$ and $r$, the function $f_0$ satisfies \autoref{eq:f0}. On the other hand,
  $g_0$ satisfies \autoref{eq:g0} because $\mathcal{B}_{L}$ is a leftward
  blend-to-zero operator of orders $\ell$ and $r$. Therefore $h_{0}$ satisfies 
  \autoref{eq:fa} and \autoref{eq:gb}  by \autoref{rmk:sumtrans}. The final assertion
  follows from \autoref{rmk:interpol}. 
\end{proof}

\begin{remark}\label{rmk:blendop}
  If $\mathcal{B}:{C}^{q}[a_{0},b_{0}]\to{C}^{q}[a_{0},b_{0}]$ is a leftward (respectively rightward) blend-to-zero operator of orders 
  $\ell$ and $r$, then $I-\mathcal{B}$ is a rightward (respectively leftward) 
  blend-to-zero operator of the same orders.
\end{remark}

\begin{corollary}\label{cor:sumblendop}
  Let $\mathcal{B}$ be a leftward blend-to-zero operator
  on $C^{q}[a_{0},b_{0}]$ of orders $\ell$ and $r$.
  Let $f\in C^{q}[a,b_{0}]$ and $g\in C^{q}[a_{0},b]$. Define $f_{0}=(I-\mathcal{B})\!\left(f|_{[a_0,b_0]}\right)$, 
  $g_{0}=\mathcal{B}\!\left(g|_{[a_0,b_0]}\right)$,
  and let $h_{0}=f_{0}+g_{0}$. 
  Then $h_0\in C^{q}[a_{0},b_{0}]$ and satisfies \autoref{eq:fa} and \autoref{eq:gb}.
  Furthermore, the piecewise function $h$ defined by \autoref{eq:pieceh} is a smooth transition function of orders $\ell$ and $r$ from $f$ to $g$ over $[a_0,b_0]$.
\end{corollary}
\begin{proof}
  The statement follows from \autoref{thm:sumblendop} and \autoref{rmk:blendop}. 
\end{proof}

Together with \autoref{rmk:interpol}, \autoref{cor:sumblendop} shows that it
suffices to construct leftward blend-to-zero operators in order to generate smooth
transition functions. In the following sections we develop several such operators.

\section{Polynomial two-point Hermite interpolation}
\label{sec:hermite}

A function $h_0$ satisfying the endpoint conditions of \autoref{eq:fa} and
\autoref{eq:gb} can be generated as a two-point
Hermite interpolating polynomial. Moreover, this interpolant can be expressed as a two-point 
Taylor expansion \cite{Estes1966,Lopez2002,Shustov2015} in which each term is 
corrected by the regularized incomplete Beta-function given by
  \begin{equation}\label{eq:Beta}
B_{\ell,r}(x) := \dfrac{\int_{0}^{x}t^\ell(1-t)^{r} dt}{\int_{0}^{1}t^\ell(1-t)^{r} dt},
\quad \ell,r\in\mathbb N_0, \quad 0\leq x\leq 1.
  \end{equation}

\begin{lemma}\label{lem:Beta}\cite[p.~305-307]{Laksa2022}
  The regularized incomplete Beta-function given by \autoref{eq:Beta} 
  is a strictly increasing function on $[0,1]$ such that:
  \begin{alignat}{7}
    B_{\ell,r}(0) & = & 0, && \quad B_{\ell,r}^{(j)}(0) & = & 0, 
    && \quad j=1,\dots,\ell, \label{eq:B0} \\
    B_{\ell,r}(1) & = & 1, && \quad B_{\ell,r}^{(k)}(1) & = & 0, 
    && \quad k=1,\dots,r. \label{eq:B1}
  \end{alignat}
\end{lemma}

\begin{figure}[htbp]
    \centering
    \begin{subfigure}[t]{0.3\textwidth}
      \centering
      \includegraphics[width=1.0\textwidth]{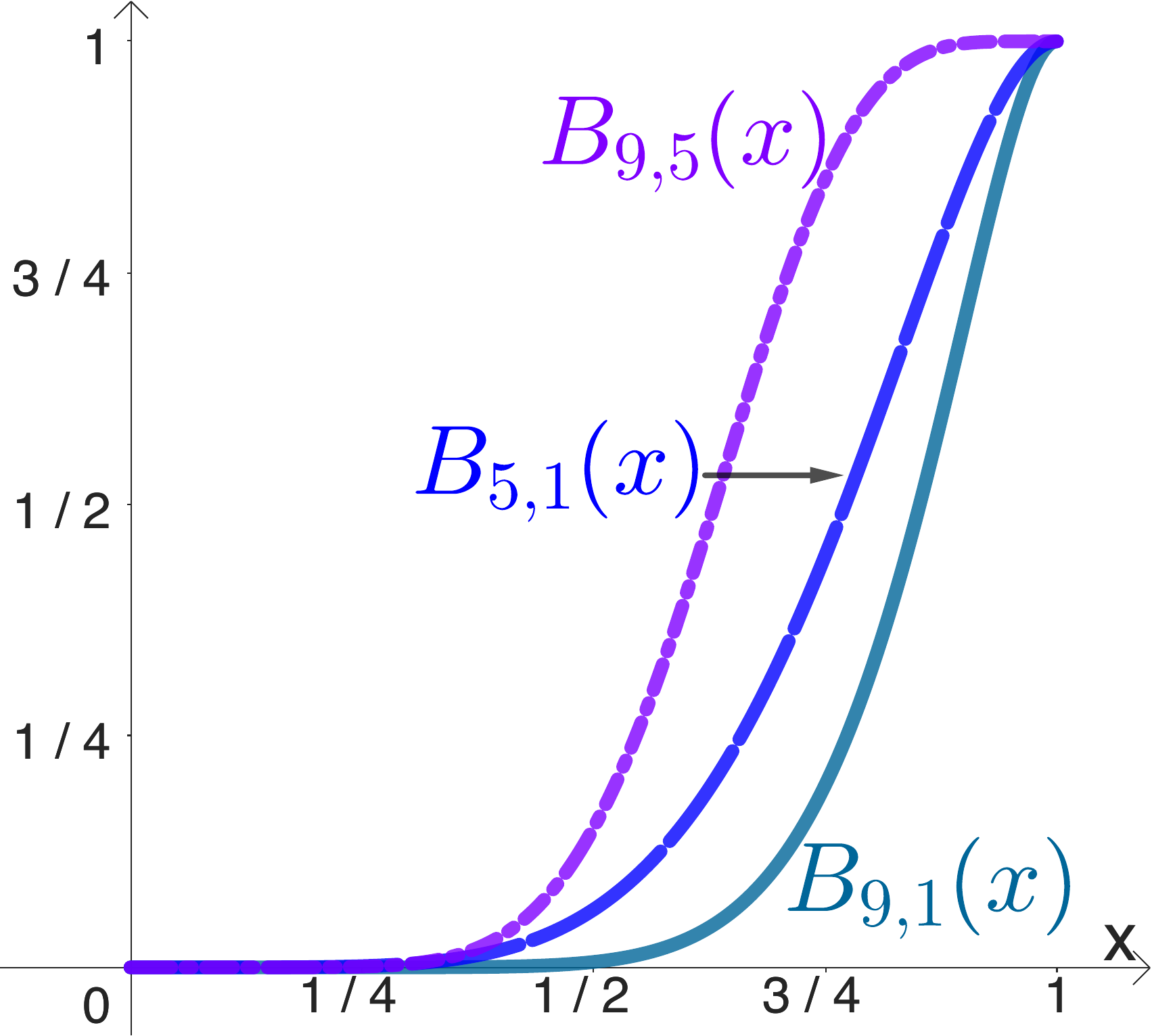}
      \caption{$B_{\ell,r}$ with $\ell>r$.}
    \end{subfigure}
    \begin{subfigure}[t]{0.3\textwidth}
      \centering
      \includegraphics[width=1.0\textwidth]{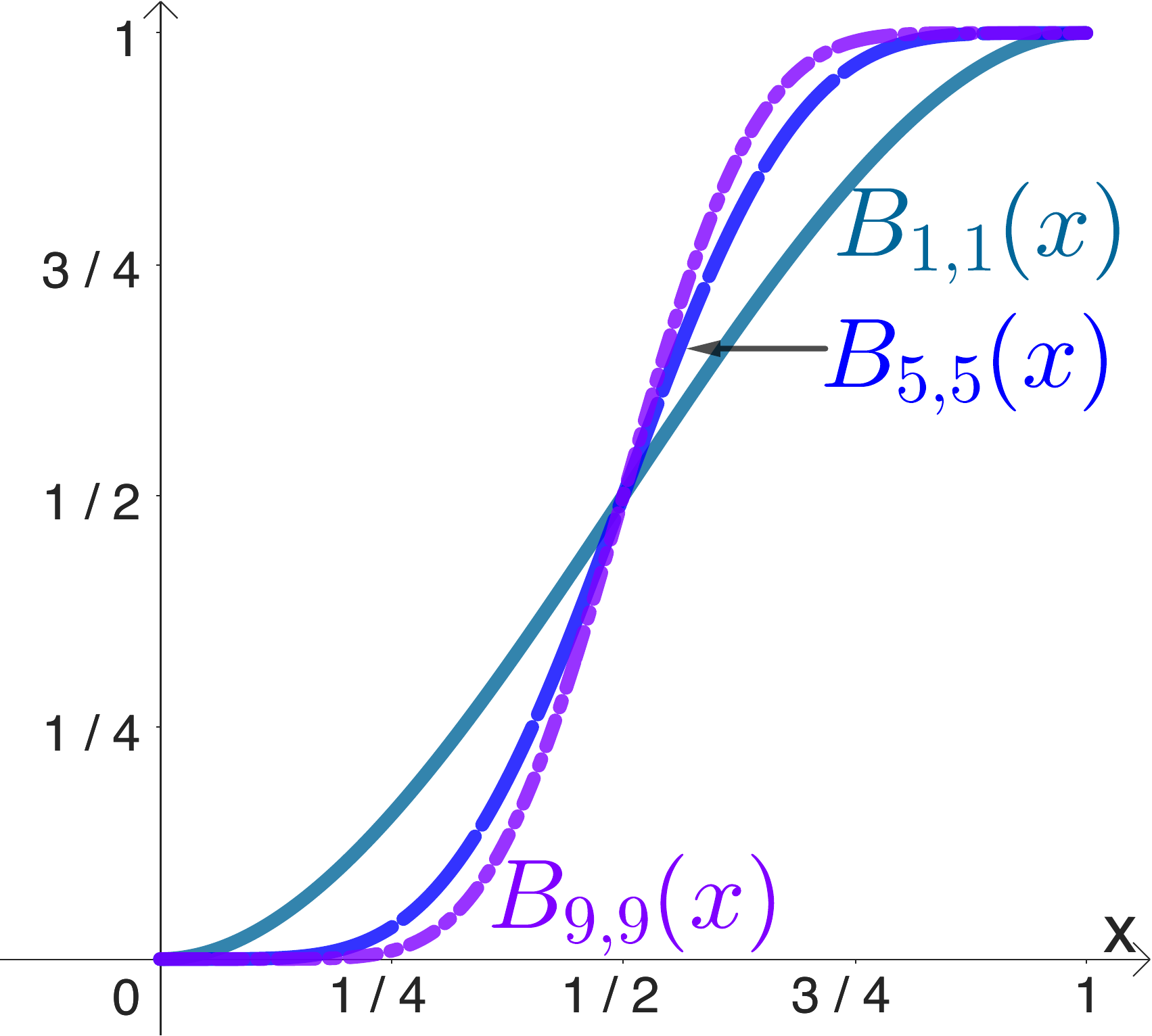}
      \caption{Symmetric $B_{\ell,r}$ ($\ell=r$).}
    \end{subfigure}
    \begin{subfigure}[t]{0.3\textwidth}
      \centering
      \includegraphics[width=1.0\textwidth]{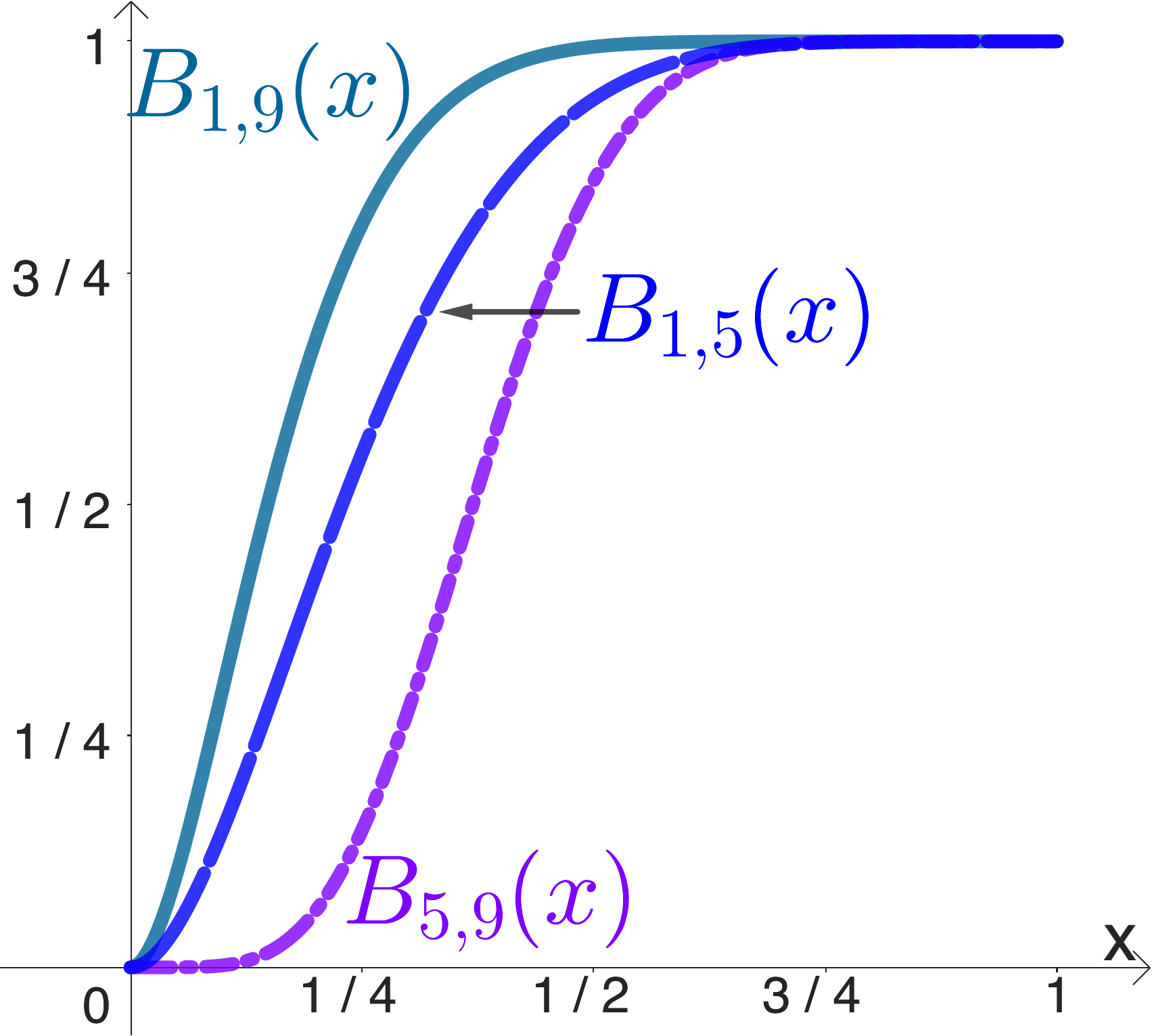}
      \caption{$B_{\ell,r}$ with $\ell<r$.}
    \end{subfigure}
    \caption{The regularized incomplete Beta-function $B_{\ell,r}$ as a 
    smooth step function of orders $\ell$ and $r$.}
    \label{fig:ibetafun}
\end{figure}

\begin{theorem}\label{thm:2Hermite}
  Let $a_{0},b_{0}\in\mathbb R$ be such that $a_{0}<b_{0}$, and let $\ell$ and $r$ be 
  non-negative integers. Let $f$ be a real-valued ${C}^{\ell}$-function 
  on a neighborhood of $x=a_{0}$, and let $g$ be a real-valued 
  ${C}^{r}$-function on a neighborhood of $x=b_{0}$. 
  Let $H$ be the polynomial of degree $\ell+r+1$ that satisfies 
  \begin{alignat}{3}
    H^{(j)}(a_{0}) & = & f^{(j)}(a_{0}), && \quad j=0,\dots,\ell, \label{eq:Ha}\\
    H^{(k)}(b_{0}) & = & g^{(k)}(b_{0}), && \quad k=0,\dots,r. \label{eq:Hb}
  \end{alignat}
  Then
  \begin{equation}\label{eq:hinterpol1}
    H(x) = 
    \sum_{j=0}^{\ell} \dfrac{f^{(j)}(a_{0})}{j!} (x-a_{0})^{j} B_{r,\ell-j}(1-\xi(x)) 
    + \sum_{k=0}^{r} \dfrac{g^{(k)}(b_{0})}{k!} (x-b_{0})^{k} B_{\ell,r-k}(\xi(x)),
  \end{equation}
 where $\xi(x)=(x-a_{0})/(b_{0}-a_{0})$. Furthermore, the regularized incomplete Beta-function can be written as 
 \begin{equation}\label{eq:Beta1}
  B_{\ell,r}(x) = x^{\ell+1}\sum_{i=0}^{r} \binom{\ell+i}{i}(1-x)^{i}.
 \end{equation}
\end{theorem}

\begin{proof}
  The Hermite interpolating polynomial of degree $\ell+r+1$ that satisfies
  \autoref{eq:Ha} and \autoref{eq:Hb} can be written as \cite{Shustov2015}
  \begin{equation}
    \begin{aligned}
      H(x) & = &  
      \left(\dfrac{b_{0}-x}{b_{0}-a_{0}}\right)^{r+1}
      \sum_{j=0}^{\ell}\dfrac{f^{(j)}(a_{0})}{j!} (x-a_{0})^{j} 
      \sum_{i=0}^{\ell-j}\binom{r+i}{i}
      \left(\dfrac{x-a_{0}}{b_{0}-a_{0}}\right)^{i}
      + \\
      & & 
      \left(\dfrac{x-a_{0}}{b_{0}-a_{0}}\right)^{\ell+1}
      \sum_{k=0}^{r}\dfrac{g^{(k)}(b_{0})}{k!} (x-b_{0})^{k}
      \sum_{i=0}^{r-k}\binom{\ell+i}{i}
      \left(\dfrac{b_{0}-x}{b_{0}-a_{0}}\right)^{i},
    \end{aligned}
  \end{equation}
  which in turn is equivalent to
  \begin{equation}\label{eq:hinterpol2}
    H(x) = \sum_{j=0}^{\ell}\dfrac{f^{(j)}(a_{0})}{j!} 
    (x-a_{0})^{j} \varphi_{r,\ell-j}(1-\xi(x)) +
    \sum_{k=0}^{r}\dfrac{g^{(k)}(b_{0})}{k!} 
    (x-b_{0})^{k}\varphi_{\ell,r-k}(\xi(x)),
  \end{equation}
  where 
  \begin{equation}\label{eq:Beta2}
    \varphi_{\ell,r}(x) = x^{\ell+1}\sum_{i=0}^{r} \binom{\ell+i}{i}(1-x)^{i}.
  \end{equation}
  In particular, substituting
  \begin{alignat}{8}
    a_{0} & = 0, \quad & f(a_{0}) & = 0, & \quad f^{(j)}(a_0) & = 0, 
    & \quad j & = 1,\dots,\ell, \\
    b_{0} & = 1, \quad & g(b_{0}) & = 1, & \quad g^{(k)}(b_0) & = 0, 
    & \quad k & = 1,\dots, r
  \end{alignat}
  into \autoref{eq:hinterpol2} yields $H=\varphi_{\ell,r}$. Then the polynomial 
  $\varphi_{\ell,r}$ satisfies \autoref{eq:B0} and \autoref{eq:B1} (replacing $B_{\ell,r}$ by $\varphi_{\ell,r}$).

  On the other hand, from  \autoref{lem:Beta} and using the {antiderivative} 
  \cite[p.~68]{Gradshteyn2007} given by
  \begin{equation}
\int t^\ell(1-t)^r dt = 
\sum_{i=0}^r \dfrac{\ell!r!}{(r-i)!(\ell+i+1)!}t^{\ell+i+1}(1-t)^{r-i},
  \end{equation}
  we have that $B_{\ell,r}$ is also a polynomial of degree 
  $\ell+r+1$ that satisfies \autoref{eq:B0} and \autoref{eq:B1}. Then from the uniqueness of the 
  Hermite interpolating polynomial, we deduce that $B_{\ell,r}=\varphi_{\ell,r}$. 
  Hence, \autoref{eq:Beta1} follows from \autoref{eq:Beta2}, and consequently, 
  \autoref{eq:hinterpol1} follows from \autoref{eq:hinterpol2}. 
\end{proof}

Linear blend-to-zero operators can be generated by two-point Hermite interpolating
polynomials.

\begin{corollary}\label{cor:hermiteblend}
  Let $a_{0},b_{0}\in\mathbb R$ be such that $a_{0}<b_{0}$, and let 
  $\ell,r,q\in\overline{\mathbb{N}}_0$ be such that $\ell,r\leq q$ and
  $\ell,r < \infty$. Let
  \begin{alignat}{5}
    \mathcal{B}_{L}(g)(x) & = & 
    \sum_{k=0}^{r} \dfrac{g^{(k)}(b_{0})}{k!}
    (x-b_{0})^{k} B_{\ell,r-k}(\xi(x)), \quad\ \
    & \quad \forall g\in {C}^{q}[a_{0},b_{0}], & \quad a_{0}\leq x\leq b_{0}, 
    \label{eq:BLhermite} \\
    \mathcal{B}_{R}(f)(x) & = & 
    \sum_{j=0}^{\ell} \dfrac{f^{(j)}(a_{0})}{j!}
    (x-a_{0})^{j} B_{r,\ell-j}(1-\xi(x)),
    & \quad \forall f\in {C}^{q}[a_{0},b_{0}], & \quad a_{0}\leq x\leq b_{0}, 
    \label{eq:BRhermite}
  \end{alignat}
  where $\xi(x)=(x-a_{0})/(b_{0}-a_{0})$. Then
  $\mathcal{B}_{L}$ is a linear leftward blend-to-zero operator
  from ${C}^{q}[a_{0},b_{0}]$ to itself of orders $\ell$ and $r$
  and $\mathcal{B}_{R}$ is a linear rightward blend-to-zero operator 
  from ${C}^{q}[a_{0},b_{0}]$ to itself of the same orders.
\end{corollary}

\begin{proof}
  Let $g\in{C}^{q}[a_{0},b_{0}]$, and let $g_{0} = \mathcal{B}_{L}(g)$. 
  Since each $B_{\ell,r-k}$ is a polynomial on $[0,1]$ and $\xi$ is affine,
we have that $B_{\ell,r-k}\circ\xi\in C^q[a_0,b_0]$ for $k=0,\dots,r$.
As $C^q[a_0,b_0]$ is a linear space closed under multiplication, each term in the sum of \autoref{eq:BLhermite}
belongs to $C^q[a_0,b_0]$, and therefore $g_0\in C^q[a_0,b_0]$.
  Thus, $\mathcal{B}_{L}$ 
  is an operator from ${C}^{q}[a_{0},b_{0}]$ to itself. Setting 
  $f^{(j)}(a_{0})=0$ for $j=0,\dots,\ell$ and $H=g_{0}$ in \autoref{thm:2Hermite}, 
  it follows that $g_{0}$ satisfies \autoref{eq:g0}. Hence, $\mathcal{B}_{L}$ is a 
  leftward blend-to-zero operator of orders $\ell$ and $r$.

  Let $g_1,g_2\in{C}^{q}[a_{0},b_{0}]$, and let $\alpha,\beta\in\mathbb R$. 
  Then for all $x\in[a_{0},b_{0}]$, we have that
  \begin{alignat*}{2} 
    \mathcal{B}_{L}(\alpha g_{1}+\beta g_{2})(x) 
    & =  
    \sum_{k=0}^{r} 
    \dfrac{(\alpha g_{1}+\beta g_{2})^{(k)}(b_{0})}{k!} 
    (x-b_{0})^{k} B_{\ell,r-k}(\xi(x)) \\ 
    & = 
    \alpha \sum_{k=0}^{r} 
    \dfrac{g_{1}^{(k)}(b_{0})}{k!} (x-b_{0})^{k} B_{\ell,r-k}(\xi(x)) +
    \beta \sum_{k=0}^{r} 
    \dfrac{g_{2}^{(k)}(b_{0})}{k!} (x-b_{0})^{k} B_{\ell,r-k}(\xi(x)) \\
    & = 
    \alpha\mathcal{B}_{L}(g_{1})(x) + \beta\mathcal{B}_{L}(g_{2})(x).
  \end{alignat*}
  Therefore, $\mathcal{B}_{L}$ is a linear operator.  The proof for the rightward operator $\mathcal{B}_{R}$ is analagous.
\end{proof}

\begin{remark}

  The degree for the polynomial in \autoref{eq:hinterpol1} increases with
the number of endpoint derivatives to be matched. Consequently, polynomial
two-point Hermite interpolation cannot generate smooth transition functions of
arbitrary order over $[a_0,b_0]$.
\end{remark}

\section{Functions with flat ends and smooth staircases}
\label{sec:staircases}
 
\begin{figure}[htbp]
	\centering
	\includegraphics[width=.5\textwidth]{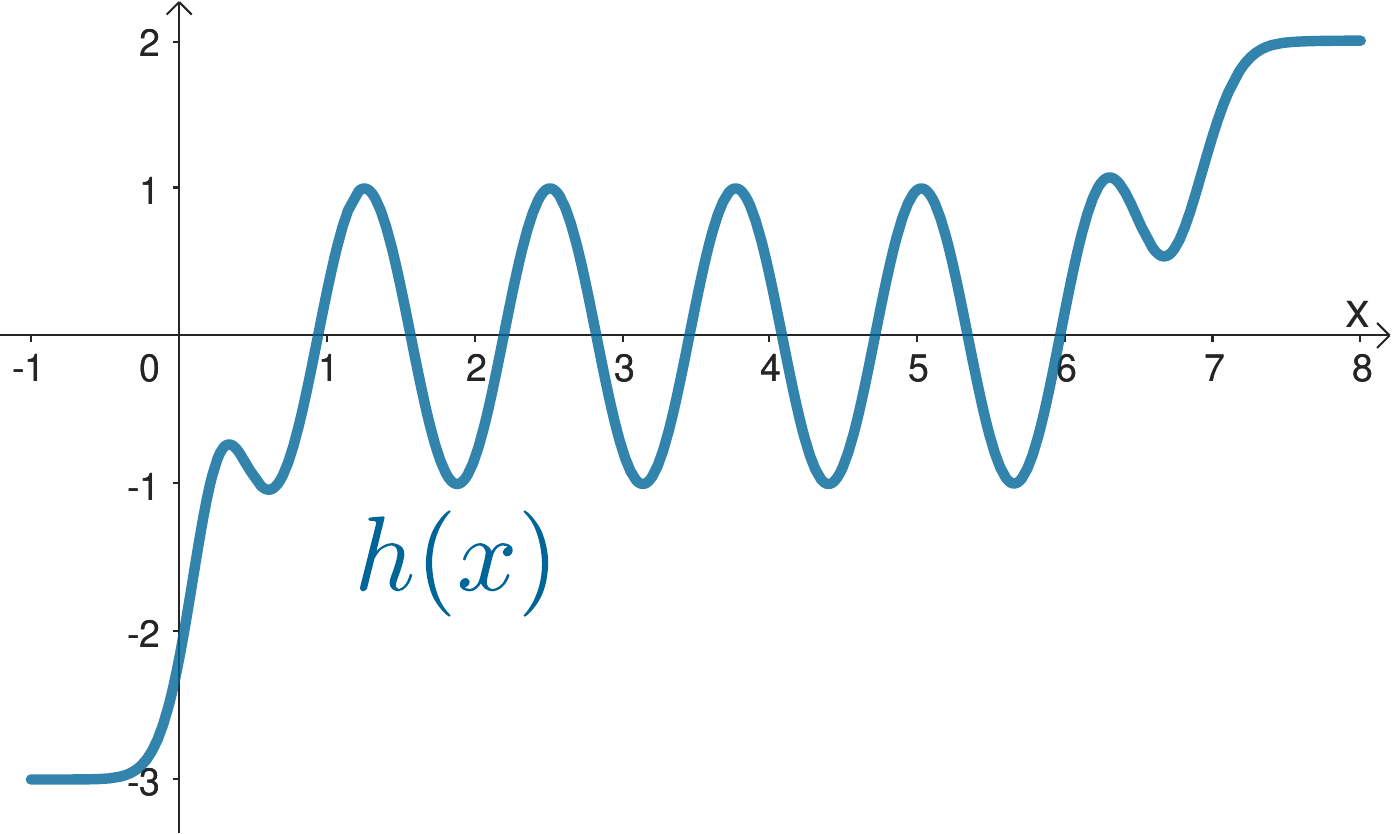}
	\caption{Illustration of a function $h\in{C}^{4}[-1,8]$ with flat ends of orders 1 and 2
  given by: $h(x) = R_{1,4}((x+1)/2)(3+\cos(5x)) - 3$ for $|x|\leq 1$; 
  $h(x)=\cos(5x)$ for $1<x<5$; and
  $h(x) = R_{4,2}((x-5)/3)(2-\cos(5x)) + \cos(5x)$
  for $5\leq x\leq 8$, where both
  $R_{1,4}$ and $R_{4,2}$ are rational B-functions \cite{Hartmann2001,Laksa2022} (presented later in \autoref{eq:rationalB}).}
	\label{fig:flatends}
\end{figure}

In order to facilitate the construction in \autoref{sec:blendop} of other blend-to-zero operators with useful features, we first examine some algebraic structures and geometric properties of \emph{functions with flat ends} on closed intervals (an interesting example of which is illustrated in \autoref{fig:flatends}.)

\begin{definition}
  Let $\ell,n,r\in\overline{\mathbb{N}}_0$ be such that $1\leq \ell,r\leq n$. A function 
  $h\in{C}^{n}[a,b]$ such that $h^{(j)}(a)=0$ for $j=1,\dots,\ell$ and 
  $h^{(k)}(b)=0$ for $k=1,\dots,r$ is called a {\em function with flat ends of 
  orders $\ell$ and $r$}. The function $h$ has order $r$ if $\ell=r$.
  In particular, $h$ has arbitrary order if $\ell=r=\infty$.
  Denote by $\mathscr{F}_{\ell,r}^{n}[a,b]$ the family of ${C}^{n}$-functions
  on $[a,b]$ with flat ends of orders $\ell$ and $r$. If $\ell=r$, the set 
  $\mathscr{F}_{\ell,r}^{n}[a,b]$ is denoted by $\mathscr{F}_{r}^{n}[a,b]$.
\end{definition}

\begin{remark}\label{rmk:FnCn}
  The set $\mathscr{F}_{\ell,r}^{n}[a,b]$ is non-empty since the constant functions 
  have flat ends of arbitrary orders. Furthermore, the following holds true:
  \begin{enumerate}[(i)]
    \item $\mathscr{F}_{\ell,r}^{n}[a,b]\subset{C}^{n}[a,b]$;
    \item $\mathscr{F}_{\ell+1,r}^{n}[a,b]\subset\mathscr{F}_{\ell,r}^{n}[a,b]$ for $\ell<n$;
    \item $\mathscr{F}_{\ell,r+1}^{n}[a,b]\subset\mathscr{F}_{\ell,r}^{n}[a,b]$ for $r<n$.
    
  \end{enumerate}
\end{remark}

Multiplication of a function by a function with flat ends preserves endpoint derivatives 
up to multiplication by a factor. 

\begin{lemma}\label{lemma:prodB}
  Let $h\in\mathscr{F}_{\ell,r}^{n}[a,b]$, and let $f\in{C}^{n}[a,b]$. 
  Then the product $hf$ belongs to ${C}^{n}[a,b]$ with 
  $(hf)^{(j)}(a)=h(a)f^{(j)}(a)$ for $j=0,\dots,\ell$
  and $(hf)^{(k)}(b)=h(b)f^{(k)}(b)$ for $k=0,\dots,r$.
\end{lemma}

\begin{proof}
  Set $\varphi:=hf$. 
  As ${C}^{n}[a,b]$ is closed under function multiplication and contains
  $\mathscr{F}_{\ell,r}^{n}[a,b]$, we deduce that $\varphi\in{C}^{n}[a,b]$. 
  Since $h\in\mathscr{F}_{\ell,r}^{n}[a,b]$, we have $h^{(i)}(a)=0$ for $i=1,\dots,\ell$ 
  and $h^{(i)}(b)=0$ for $i=1,\dots,r$. Then from the general Leibniz rule given by
  $\varphi^{(m)} = \sum_{i=0}^m \binom{m}{i}h^{(i)}f^{(m-i)}$
  for $m=0,\dots,n$, it follows that
  $\varphi^{(j)}(a)=h(a)f^{(j)}(a)$ for $j=0,\dots,\ell$ and 
  $\varphi^{(k)}(b)=h(b)f^{(k)}(b)$ for $k=0,\dots,r$.  
\end{proof}

Functions with flat ends admit useful algebraic structures and properties as follows.

\begin{theorem}\label{thm:prodB}
  The set $\mathscr{F}_{\ell,r}^{n}[a,b]$ is closed under function multiplication.
\end{theorem}

\begin{proof}
  Let $h,f\in\mathscr{F}_{\ell,r}^{n}[a,b]$, let $j\in\{1,\dots,\ell\}$, 
  let $k\in\{1,\dots,r\}$, and set $\varphi:=hf$. Then by \autoref{lemma:prodB}, 
  we have that $\varphi\in{C}^{n}[a,b]$ such that 
  $\varphi^{(j)}(a)=h(a)f^{(j)}(a)$ and $\varphi^{(k)}(b)=h(b)f^{(k)}(b)$. 
  Since $f\in\mathscr{F}_{\ell,r}^{n}[a,b]$, we have that $f^{(j)}(a)=0$ and 
  $f^{(k)}(b)=0$. Consequently, $\varphi^{(j)}(a)=0$ and $\varphi^{(k)}(b)=0$. 
  Therefore, $\varphi\in\mathscr{F}_{\ell,r}^{n}[a,b]$. 
\end{proof}

\begin{theorem}\label{thm:lblend}
  The set $\mathscr{F}_{\ell,r}^{n}[a,b]$ is a linear subspace of ${C}^{n}[a,b]$.
\end{theorem}

\begin{proof}
  Let $f,g\in\mathscr{F}_{\ell,r}^{n}[a,b]$, let $\alpha_1,\alpha_2\in\mathbb R$, 
  and set $h:=\alpha_{1}f+\alpha_{2}g$. Then $h\in{C}^{n}[a,b]$ 
  because ${C}^{n}[a,b]$ is a linear space that contains $\mathscr{F}_{\ell,r}^{n}[a,b]$. 
  Now let $j\in\{1,\dots,\ell\}$, and let $k\in\{1,\dots,r\}$. Since both $f$ 
  and $g$ are ${C}^{n}$-functions on $[a,b]$ with flat ends of orders  $\ell$ and $r$, we have that $f^{(j)}(a) = g^{(j)}(a) = 0$ and 
  $f^{(k)}(b) = g^{(k)}(b)=0$. Hence $h^{(j)}(a)=0$ and $h^{(k)}(b)=0$. 
  Therefore, $h\in\mathscr{F}_{\ell,r}^{n}[a,b]$. 
\end{proof}

Inspired by \cite{Laksa2022,Leopold2020,Li2004}, we define and construct smooth staircase and step 
functions as monotone functions with flat ends.

\begin{definition}
  Let $\ell,n,r\in\overline{\mathbb{N}}_0$ be such that $\ell,r\leq n$. A monotone function 
  $s:[a,b]\to\mathbb{R}$ such that $s\in\mathscr{F}_{\ell,r}^{n}[a,b]$  is called a 
  {\em ${C}^{n}$-staircase function on $[a,b]$ of orders $\ell$ and $r$}.
  A strictly increasing $C^\infty$-staircase function $\sigma: [0,1] \to [0,1]$ of orders $\ell$ and $r$ such that $\sigma(0)=0$ and $\sigma(1)=1$ is called a {\em smooth step function of orders $\ell$ and $r$}.
  A ${C}^{n}$-staircase function or a smooth step function has order $r$ 
  if $\ell=r$. Denote by $\Sigma_{\ell,r}$ the family of smooth step functions 
  of orders $\ell$ and $r$, where $\Sigma_r = \Sigma_{\ell,r}$ if $\ell=r$.
\end{definition}

\begin{remark}
Although $\sigma \in C^\infty[0,1]$, it does not necessarily admit an extension to $C^\infty(\mathbb{R})$. Indeed, if $\sigma \in \Sigma_r$, then the piecewise function given by
\begin{equation}
  \widetilde{\sigma} = \begin{cases} 0, & x \leq 0\\ \sigma(x), & 0<x<1\\ 1, & x\geq 1,\end{cases}
\end{equation}
belongs to $C^r(\mathbb{R})$.
\end{remark}

Some examples of smooth step functions known in literature are the following:
\begin{enumerate}
  \item From \autoref{lem:Beta} and \autoref{thm:2Hermite}, we have that the regularized incomplete 
  Beta-function $B_{\ell,r}$ given by \autoref{eq:Beta} is a polynomial step function 
  of orders $\ell$ and $r$, as shown in \autoref{fig:ibetafun}. It can also be considered
  as a polynomial approximation to the signum function \cite{Hosenthien1972}.
  \item The rational B-function \cite{Hartmann2001,Laksa2022}, given by
    \begin{equation}\label{eq:rationalB}
      R_{\ell,r}(x) = \dfrac{x^{\ell+1}}{x^{\ell+1}+(1-x)^{r+1}}, 
      \quad 0\leq x\leq 1,
    \end{equation}
    is a smooth step function of orders $\ell$ and $r$.
  \item The logistic expo-rational B-function \cite{Dechevsky2013}, given by
    \begin{equation}
      \mathcal{E}(x) = 
      \begin{cases}
      0, & x = 0;\\ \dfrac{1}{1+\exp\left(\frac{1}{x}-\frac{1}{1-x}\right)}, & 0<x<1; \\ 1, & x = 1,
      \end{cases}
    \end{equation}
    is a smooth step function of arbitrary orders that is not analytic at $x=0$ and $x=1$.
  \item The Fabius function $\mathcal{T}$ on $[0,1]$, defined by the equation \cite{Fabius1966} and given by
    \begin{equation}
      \mathcal{T}(x) = 
      \begin{cases}
        \int_{0}^{2x}\mathcal{T}(t)dt, & 0\leq x\leq\frac{1}{2}; \\
        \int_{2x-1}^{1}\mathcal{T}(t)dt + 2x-1, & \frac{1}{2}<x\leq 1,
      \end{cases}
    \end{equation}
    is a smooth step function of arbitrary orders that is nowhere analytic.
\end{enumerate}

The products and compositions of smooth step functions are in turn smooth step 
functions of the same orders. In order to verify this claim, we first show that 
a change of interval transforms a function with flat ends 
into another one of the same orders.

\begin{lemma}\label{lem:compB}
  Let $h\in\mathscr{F}_{\ell,r}^{n}[c,d]$, and let $f$ be a 
  ${C}^{n}$-function from $[a,b]$ to $[c,d]$ such that 
  $f(a) = c$ and $f(b) = d$. Then $h\circ f$ belongs to $\mathscr{F}_{\ell,r}^{n}[a,b]$.
\end{lemma}

\begin{proof}
  Set $\varphi:=h\circ f$. Then $\varphi\in{C}^{n}[a,b]$ because
  $f$ is a ${C}^{n}$-function from $[a,b]$ to $[c,d]$ and 
  $h\in{C}^{n}[c,d]$. Since $h$ is a ${C}^{n}$-function 
  on $[c,d]$ with flat ends of orders $\ell$ and $r$, we have that 
$f(a) = c$ and $f(b) = d$ implies $(h^{(p)}\circ f)(a)=0$ for $p=1,\dots,\ell$ 
  and $(h^{(p)}\circ f)(b)=0$ for $p=1,\dots,r$.
  Then from the general chain rule \cite[p.~22]{Gradshteyn2007}, 
  $\varphi^{(m)}=\sum_{p=1}^m U_{m,p}/p! \cdot h^{(p)}\circ f$ where
  $U_{m,p} = \sum_{i=0}^{p-1} (-1)^{i}\binom{p}{i}f^{i}\cdot(f^{p-i})^{(m)}$
  for $m=1,\dots,n$, we deduce that $\varphi^{(j)}(a)=0$ for $j=1,\dots,\ell$ 
  and $\varphi^{(k)}(b)=0$ for $k=1,\dots,r$. Therefore,
  $\varphi\in\mathscr{F}_{\ell,r}^{n}[a,b]$. 
\end{proof}

\begin{theorem}
  The set $\Sigma_{\ell,r}$ is closed under both multiplication 
  and composition of functions.
\end{theorem}

\begin{proof}
  Let $h,f\in\Sigma_{\ell,r}$. Then, by definition, both $h$ and $f$ are strictly increasing 
  ${C}^{\infty}$-functions from $[0,1]$ to $[0,1]$ with flat ends 
  of orders $\ell$ and $r$. Consequently, by \autoref{thm:prodB}, we have that 
  $hf\in\mathscr{F}_{\ell,r}^{\infty}[0,1]$. Similarly, 
  \autoref{lem:compB} implies that $h\circ f\in\mathscr{F}_{\ell,r}^{\infty}[0,1]$. Furthermore, since $h(0)=f(0)=0$ and $h(1)=f(1)=1$, then both $hf$ and $h\circ f$
map $0$ to $0$ and $1$ to $1$.
  Since both the composition and multiplication of any two strictly increasing functions 
  from $[0,1]$ to $[0,1]$ are in turn strictly increasing functions from $[0,1]$ to 
  $[0,1]$, we conclude that $hf\in\Sigma_{\ell,r}$ and $h\circ f\in\Sigma_{\ell,r}$. 
\end{proof}

\begin{remark}\label{rmk:affine}
  Let $h\in\mathscr{F}_{\ell,r}^{n}[a,b]$. The reflection over the $x$-axis flips 
  the sign of $h$, a vertical translation adds a constant to $h$, and a vertical 
  stretching multiplies $h$ by a positive constant. Then, from \autoref{thm:lblend}, 
  it follows that these operations transform $h$ into another function in 
  $\mathscr{F}_{\ell,r}^{n}[a,b]$. On the other hand, the reflection of $h$ over 
  the $y$-axis is in $\mathscr{F}_{r,\ell}^{n}[-b,-a]$, the horizontal stretching 
  of $h$ by $\alpha>0$ is in $\mathscr{F}_{\ell,r}^{n}[a/\alpha, b/\alpha]$, 
  and the horizontal translation of $h$ by $x_{0}\in\mathbb{R}$ is in 
  $\mathscr{F}_{\ell,r}^{n}[a+x_0,b+x_0]$. 
  
  Indeed, since both stretching and translation preserve monotonicity while axis-aligned reflections reverse it, these operations transform ${C}^n$-staircase functions into others 
  of the same orders on the corresponding intervals (with the endpoint orders interchanged in the case of reflection over the $y$-axis, as mentioned above).
\end{remark}

Smooth step functions can be extended to smooth staircase functions 
on any closed interval.

\begin{theorem}\label{cor:incblend}
  Let $a,b,c,d\in\mathbb R$ be such that $a<b$ and $c<d$.
  Let $\sigma\in\Sigma_{\ell,r}$, and set 
  \begin{equation}s(x):= c + (d-c)\sigma((x-a)/(b-a)) \qquad a\leq x\leq b. 
  \end{equation}
  Then $s$ is an increasing ${C}^{\infty}$-staircase function 
  from $[a,b]$ to $[c,d]$ of orders $\ell$ and $r$. 
\end{theorem}

\begin{proof}
  Set $\lambda(x):=(x-a)/(b-a)$ for $x\in[a,b]$ and $\tau(x):=c+(d-c)x$ for 
  $x\in[0,1]$. Given that $\lambda$ maps $[a,b]$ onto $[0,1]$, $\sigma$ maps 
  $[0,1]$ onto itself, and $\tau$ maps $[0,1]$ onto $[c,d]$, we deduce that 
  $s=\tau\circ\sigma\circ\lambda$ and $s$ maps $[a,b]$ onto $[c,d]$.

  Since $\sigma\in\Sigma_{\ell,r}$, we have that
  $\sigma\in\mathscr{F}_{\ell,r}^{\infty}[0,1]$, 
  and since $\lambda$ is a smooth function from 
  $[a,b]$ onto $[0,1]$ such that $\lambda(a)=0$ and $\lambda(b)=1$,
  \autoref{lem:compB} implies
  $\sigma\circ\lambda\in\mathscr{F}_{\ell,r}^{\infty}[a,b]$. 
  As $\tau$ stretches and translates $\sigma\circ\lambda$, we deduce from 
  \autoref{rmk:affine} that $s\in\mathscr{F}_{\ell,r}^{\infty}[a,b]$.
  Noting that $\tau$, $\sigma$ and $\lambda$ are increasing functions
  and using the fact that the composition of increasing functions is in turn 
  an increasing function, we deduce that $s$ is an increasing function (and hence monotone from $[a,b]$ onto $[c,d]$). Together with
  $s\in\mathscr{F}_{\ell,r}^{\infty}[a,b]$, this proves the assertion.
\end{proof}

Motivated by \cite{Laksa2022}, we next generate smooth step and staircase functions
with point-symmetry.

\begin{definition}
     A function $u:[a,b]\to\mathbb{R}$ is symmetric on $[a,b]$ if its graph is
  symmetric about the point $(x_M,u(x_M))$, where $x_M=(a+b)/2$.
\end{definition}

\begin{figure}[htbp]
    \centering
    \begin{subfigure}[t]{0.45\textwidth}
      \centering
      \includegraphics[width=0.7\textwidth]{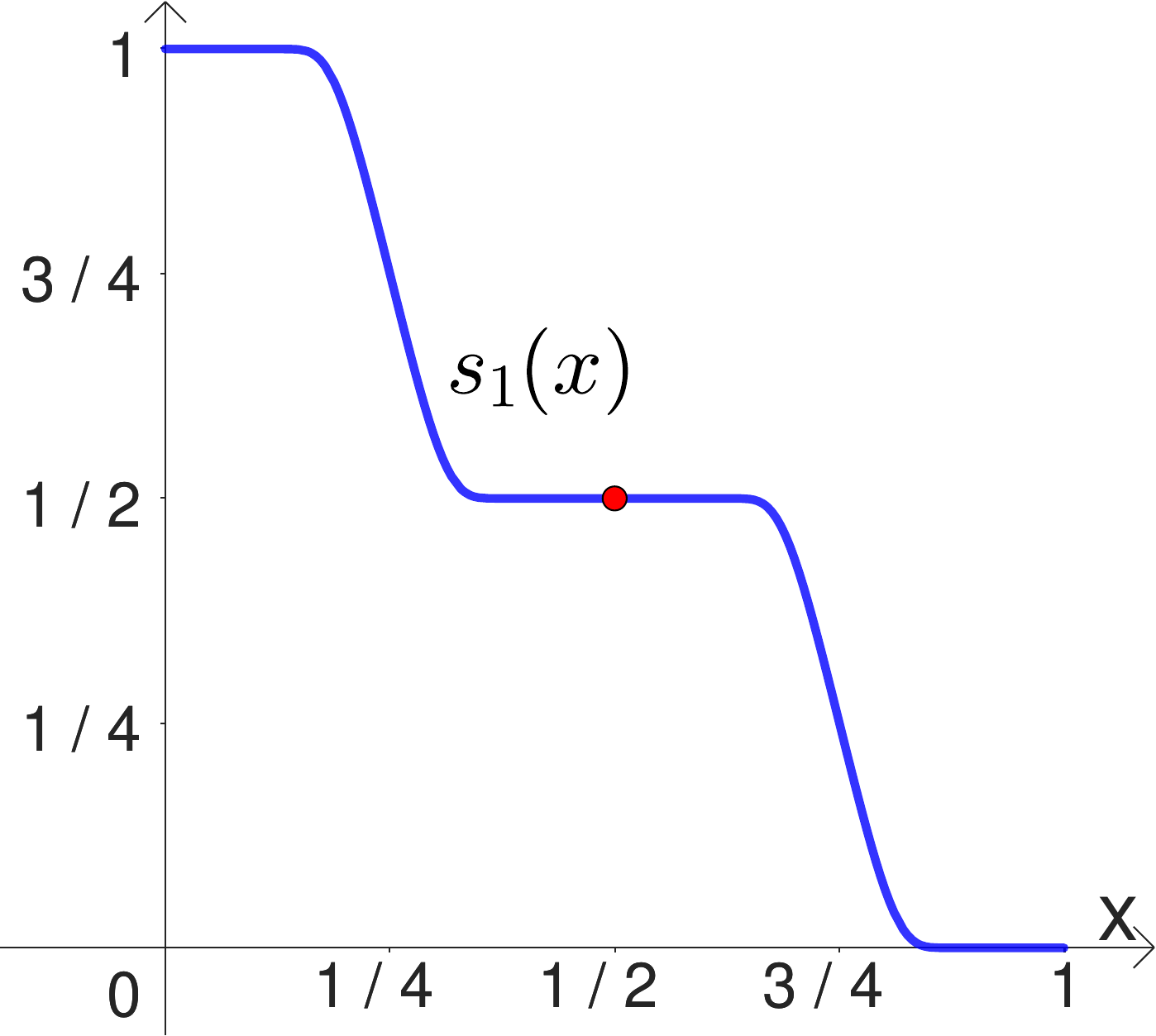}
      \caption{${C}^{\infty}$-staircase function of arbitrary order.}
    \end{subfigure}
    \begin{subfigure}[t]{0.45\textwidth}
      \centering
      \includegraphics[width=0.7\textwidth]{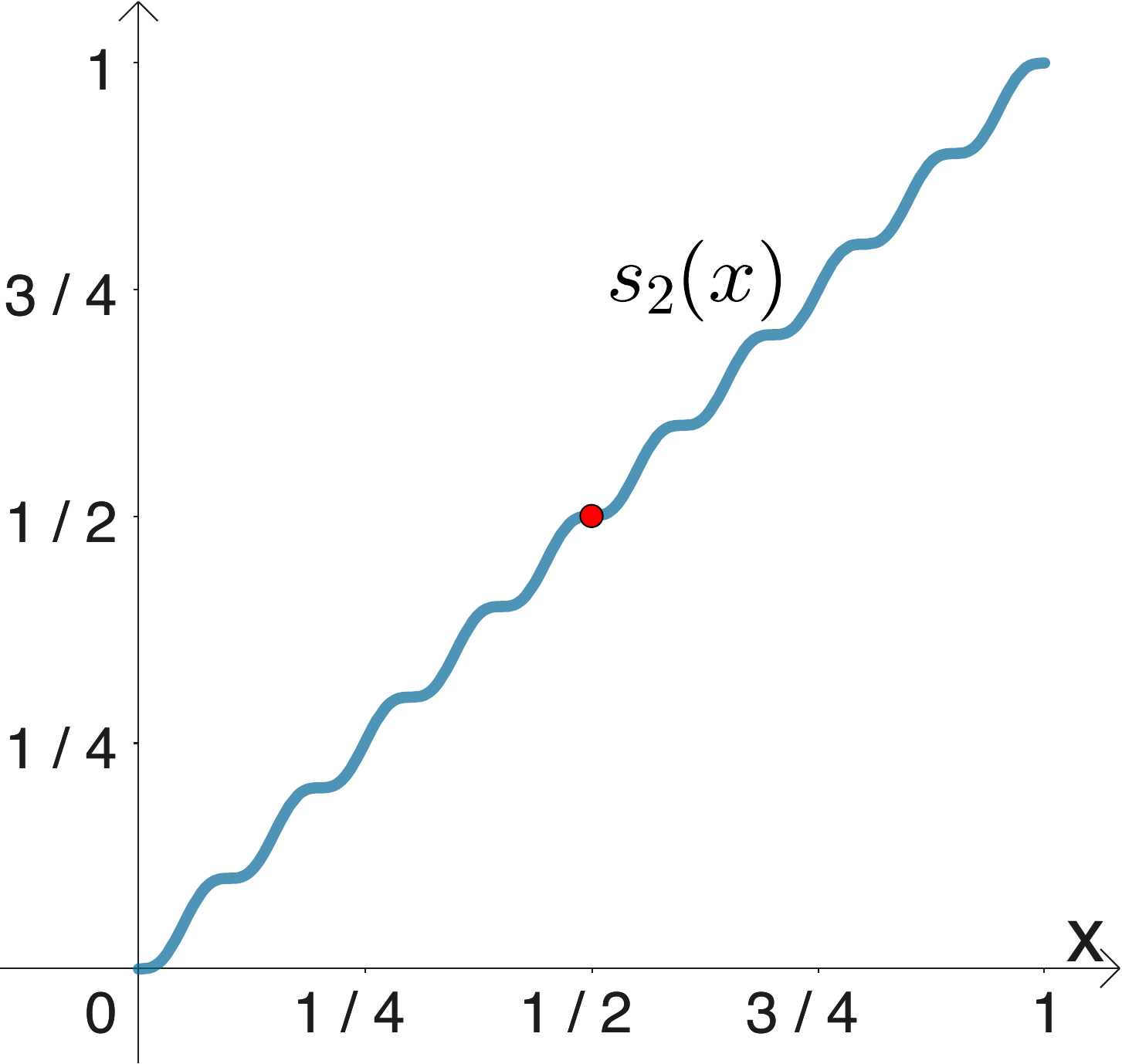}
      \caption{${C}^{2}$-staircase function of order 2.}
    \end{subfigure}
    \caption{Illustrative examples of staircase functions on $[0,1]$ symmetric about the red midpoint: (a) $s_{1}(x) = 1-\mathcal{E}(\mathcal{E}(2x))/2$ for $0\leq x< 1/2$
    and  $s_{1}(x) = 1/2-\mathcal{E}(\mathcal{E}(2x-1))/2$ for $1/2\leq x\leq 1$,
    where $\mathcal{E}$ is the logistic expo-rational B-function; and 
    (b) $s_{2}(x)=x-\sin(20\pi x)/(20\pi)$.}
    \label{fig:staircase}
\end{figure}

Examples of symmetric staircase functions are presented in \autoref{fig:staircase}.

\begin{lemma}\label{lemma:sym}
  Let $u$ be a real-valued function on $[a,b]$. Then $u$ is symmetric on $[a,b]$ 
  if and only if $u(x) + u(a+b-x) = u(a)+u(b)$ for $a\leq x\leq b$.
\end{lemma}

\begin{proof}
 Let $x_{M}=(a+b)/2$, and let $y_{M}=u(x_{M})$.
 Since the reflection $r$ of $u$ about the point $(x_{M},y_{M})$
 is given by $r(x)=2y_{M}-u(2x_{M}-x)$ for $a\leq x\leq b$,
 we have that $u$ is symmetric on $[a,b]$ if and only if
 $u(x)=2y_{M}-u(a+b-x)$ for $a\leq x\leq b$.
 By substituting either $x=x_{M}$ in $u(x) + u(a+b-x) = u(a)+u(b)$ or $x=a$ in $u(x)=2y_{M}-u(a+b-x)$, we deduce that
 $2y_M = u(a)+u(b)$. Hence, the assertion holds true. 
\end{proof}

A symmetric function on a closed interval that is flat at one end must be flat 
at the other end. Furthermore, a strictly increasing symmetric function on the 
unit interval that is flat at the origin is a smooth step function.

\begin{theorem}\label{thm:flatsym}
  Let $\ell,n\in\overline{\mathbb{N}}_0$ be such that $1\leq\ell\leq n$.
  Let $u$ be a symmetric function on $[a,b]$ such that $u\in C^{n}[a,b]$ 
  and $u^{(k)}(a)=0$ for $k=1,\dots,\ell$. Then $u\in\mathscr{F}_{\ell}^{n}[a,b]$. 
\end{theorem}

\begin{proof}
  Let $k\in\{1,\dots,\ell\}$, and let $x\in[a,b]$. Using \autoref{lemma:sym} and the fact that 
  $u$ is symmetric on $[a,b]$, we have that $u(x)+u(a+b-x)=u(a)+u(b)$. Then by induction 
  and the chain rule we deduce that $u^{(k)}(x)=(-1)^{k+1}u^{(k)}(a+b-x)$. Consequently, 
  $u^{(k)}(a)=0$ implies $u^{(k)}(b)=0$. Hence, $u$ is a ${C}^{n}$-function on $[a,b]$ with flat ends 
  of order $\ell$. 
\end{proof}

\begin{corollary}\label{cor:symstep}
  Let $\ell\in\overline{\mathbb{N}}_0$ be such that $\ell\geq1$.
  Let $\sigma\in{C}^{\infty}[0,1]$ be such that:
   \begin{enumerate}[(i)]
    \item $\sigma(x)+\sigma(1-x)=1$ for $0\leq x\leq 1$;
    \item $\sigma^{(k)}(0)=0$ for $k=0,\dots,\ell$;
    \item $\sigma'(x)>0$ for $0<x<1$.
   \end{enumerate}
   Then $\sigma$ is a smooth step function.
\end{corollary}

\begin{proof}
  From \autoref{lemma:sym} and $(i)$, it follows that $\sigma$ is symmetric on $[0,1]$. 
  Then from \autoref{thm:flatsym} and $(ii)$ we deduce that 
  $\sigma\in\mathscr{F}_{\ell}^{\infty}[0,1]$. Since $(i)$ and $(ii)$ imply 
  $\sigma(0)=0$ and $\sigma(1)=1$, and since $(iii)$ implies that $\sigma$ 
  is a strictly increasing function on $[0,1]$, then $\sigma$ maps $[0,1]$ 
  onto itself. 
  Therefore, $\sigma$ is a smooth step function of order $\ell$. 
\end{proof}

Compositions of symmetric functions are in turn symmetric functions. In particular, 
a linear change of interval preserves the symmetry of a smooth step function.

\begin{theorem}\label{thm:comsym}
  Let $u:[c,d]\to\mathbb{R}$ be a symmetric function on $[c,d]$, and let 
  $v:[a,b]\to[c,d]$ be a symmetric function on $[a,b]$ such that $v(a)=c$ and $v(b)=d$.
  Then $u\circ v$ is a symmetric function on $[a,b]$.
\end{theorem}

\begin{proof}
  Let $x\in[a,b]$. Since $u$ is symmetric on $[c,d]$, \autoref{lemma:sym} implies that
  $u(c+d-v(x)) = u(c)+u(d)-u(v(x))$. On the other hand, since $v$ is symmetric on 
  $[a,b]$, we deduce from \autoref{lemma:sym} that $v(a+b-x) = v(a)+v(b)-v(x)$. 
  Then, from $v(a)=c$ and $v(b)=d$, it follows that 
  $u(v(a+b-x)) = u(v(a))+u(v(b))-u(v(x))$.
  From this, \autoref{lemma:sym} implies that $u\circ v$ is a symmetric function on $[a,b]$. 
\end{proof}

\begin{corollary}
  Let $a,b\in\mathbb R$ be such that $a<b$ and let $\sigma$ be a symmetric function on $[0,1]$ such that $\sigma\in\Sigma_{\ell}$. 
 Then $\sigma\circ\lambda$ is an increasing ${C}^{\infty}$-staircase 
 function from $[a,b]$ to $[0,1]$ of order $\ell$ that is symmetric on $[a,b]$, 
 where $\lambda(x)=(x-a)/(b-a)$ for $a\leq x\leq b$.
\end{corollary}

\begin{proof}
  From \autoref{lemma:sym} and the fact that $\lambda$ satisfies $\lambda(a)=0$, 
  $\lambda(b)=1$ and $\lambda(x)+\lambda(a+b-x)=1$ for $a\leq x\leq b$, 
  we have that $\lambda$ is symmetric on $[a,b]$. Then, from \autoref{thm:comsym}, 
  we deduce that $\sigma\circ \lambda$ is symmetric on $[a,b]$.
  From this, \autoref{cor:incblend} implies that $\sigma\circ\lambda$ is an increasing  
  ${C}^{\infty}$-staircase function from $[a,b]$ to $[0,1]$ 
  of order $\ell$. 
\end{proof}

\section{Blend-to-zero operators from functions with flat ends}
\label{sec:blendop}

Linear blend-to-zero operators can be generated through multiplication by functions with flat ends.

\begin{theorem}\label{thm:blendop}
  Let $\ell,m,n,r\in\overline{\mathbb{N}}_0$ be such that $1\leq \ell,r\leq n\leq m$.
  Let $h\in\mathscr{F}_{\ell,r}^{m}[a,b]$, and set $\mathcal{B}(f):=hf$ for 
  $f\in{C}^{n}[a,b]$. Then 
  \begin{enumerate}[(i)]
    \item $\mathcal{B}$ is a linear operator from 
  ${C}^{n}[a,b]$ into itself;
    \item $\mathcal{B}$ maps $\mathscr{F}_{\ell,r}^{n}[a,b]$ into itself;
    \item If $h(a)=0$ and $h(b)=1$, then $\mathcal{B}$ is a leftward 
    blend-to-zero operator of orders $\ell$ and $r$;
    \item If $h(a)=1$ and $h(b)=0$, then $\mathcal{B}$ is a rightward 
    blend-to-zero operator of orders $\ell$ and $r$.
  \end{enumerate}
\end{theorem}

\begin{proof}
  
  (i) Given that $h\in\mathscr{F}_{\ell,r}^{m}[a,b]$ with $n\leq m$, 
  then $h$ is also an element of $\mathscr{F}_{\ell,r}^{n}[a,b]$, and hence
  \autoref{rmk:FnCn} implies that $h\in{C}^{n}[a,b]$.
  As ${C}^{n}[a,b]$ is closed under function multiplication, 
  it follows that $\mathcal{B}$ maps ${C}^{n}[a,b]$ into itself. 
  The operator $\mathcal{B}$ is linear because of the distributive 
  property of both function multiplication and scalar multiplication for 
  the addition of functions in ${C}^{n}[a,b]$. Assertion {\em(ii)} follows from \autoref{thm:prodB} since $\mathcal{B}$ multiplies two 
  functions in $\mathscr{F}_{\ell,r}^{n}[a,b]$. 

  Let $f\in {C}^{n}[a,b]$, and set $f_0=\mathcal{B}(f)$. 
  From \autoref{lemma:prodB} and $h\in\mathscr{F}_{\ell,r}^{n}[a,b]$, 
  it follows that $h(a)=0$ and $h(b)=1$ imply 
  $f_0^{(j)}(a)=0$ for $j=0,\dots,\ell$ and $f_0^{(k)}(b)=f^{(k)}(b)$ for $k=0,\dots,r$,
  while $h(a)=1$ and $h(b)=0$ 
  imply $f_0^{(j)}(a)=f^{(j)}(a)$ for $j=0,\dots,\ell$ and $f_0^{(k)}(b)=0$ for $k=0,\dots,r$.
  Then, from $h(a)=0$ and $h(b)=1$, we deduce that $\mathcal{B}$ is a leftward 
  blend-to-zero operator of orders $\ell$ and $r$. Similarly, from $h(a)=1$ and $h(b)=0$, 
  we have that $\mathcal{B}$ is a rightward blend-to-zero 
  operator of the same orders. Hence, assertions {\em(ii)} and {\em(iii)} 
  hold true. 
\end{proof}

It suffices to generate smooth step functions in order to construct
linear blend-to-zero operators as shown in \autoref{fig:blend0op}.

\begin{figure}[htbp]
	\centering
	\includegraphics[width=.5\textwidth]{blend0op.pdf}
\caption{Illustration of the leftward blend-to-zero operator
  $\mathcal{B}_{L}:C^{n}[2,4]\to C^{n}[2,4]$ 
  given by $f\mapsto f_{0}$, where $f_{0}(x)=f(x)\sigma(x/2-1)$ and 
  $\sigma$ is the rational B-function $R_{4,2}$ (see~\autoref{eq:rationalB}). The operator flattens 
  $f(x)=2+(5-x)\cos^{2}(3\pi(5-x))$ at $x=2$.}
	\label{fig:blend0op}
\end{figure}

\begin{corollary}\label{cor:leftblend}
   Let $a,b\in\mathbb R$ be such that $a<b$, and let $\ell,n,r\in\overline{\mathbb{N}}_0$ be such that $\ell,r\leq n$. 
  Let $\sigma\in\Sigma_{\ell,r}$, and set 
  $\mathcal{B}_{L}(f):=(\sigma\circ\lambda)\cdot f$ for $f\in{C}^{n}[a,b]$, 
  where $\lambda(x):=(x-a)/(b-a)$ for $a\leq x\leq b$. Then $\mathcal{B}_{L}$ is a linear leftward blend-to-zero operator on ${C}^n[a,b]$ of 
  orders $\ell$ and $r$.
\end{corollary}

\begin{proof}
  Let $h=\sigma\circ\lambda$. From \autoref{cor:incblend} and $\sigma\in\Sigma_{\ell,r}$, 
  we have that $h$ is a strictly increasing function from $[a,b]$ onto $[0,1]$
  such that $h\in\mathscr{F}_{\ell,r}^{\infty}[a,b]$ and $h(a)=0$, $h(b)=1$.
  Then the assertion follows from \autoref{thm:blendop}. 
\end{proof}

\section{Trigonometric smooth step functions}
\label{sec:trigstep}

We can generate smooth step functions of odd order by integrating odd powers of sine.

\begin{theorem}\label{thm:trigstep}
 Let $m$ be a non-negative integer, let $S_m(x):= \int_{0}^{x}\sin^{2m+1}(\pi t)\,dt$ for 
 $0\leq x\leq 1$, and let $\mathcal{T}_m(x):={S_m(x)}/{S_m(1)}$ for $0\leq x\leq 1$.
 Then $\mathcal{T}_m$ is a symmetric function on $[0,1]$ such that 
 $\mathcal{T}_m\in\Sigma_{2m+1}$.
 
\end{theorem}

\begin{proof}
  We first note that $S_m(x)>0$ for $0<x\leq 1$ because of
  $\sin(\pi t)>0$ for $0<t<1$. This guarantees that 
  $\mathcal{T}_m$ is a non-negative function defined on $[0,1]$.
  Furthermore, by the Fundamental Theorem of Calculus, we have
  \begin{equation}\label{eq:dFm}
    \mathcal{T}'_m(x)=\sin^{2m+1}(\pi x)/S_m(1),\quad 0\leq x\leq 1.
  \end{equation}
Since sine is a smooth function such that $\sin(\pi x)>0$ for $0<x<1$, 
it follows from \autoref{eq:dFm} that $\mathcal{T}_m\in C^{\infty}[0,1]$
  such that $\mathcal{T}_m'(x)>0$ for $0<x<1$. Noting that 
  $\sin(\pi(1-t)) = \sin(\pi t)$ for $t\in\mathbb R$ and using integration by 
  substitution, we deduce that $S_m(1-x)=S_m(1) - S_m(x)$ for $0\leq x\leq 1$, 
  hence $\mathcal{T}_m(1-x)=1-\mathcal{T}_m(x)$ for $0\leq x\leq 1$. 
  Then, from \autoref{lemma:sym}, it follows that $\mathcal{T}_m$ is a symmetric 
  function on $[0,1]$.

Since the leading-order term in the Maclaurin series of $\sin(\pi x)$ is $\pi x$,
the leading-order term in the Maclaurin series of $\sin^{2m+1}(\pi x)$ is $\left(\pi x\right)^{2m+1}$.
  Thus, the derivatives of $\sin^{2m+1}(\pi x)$ at $x=0$ 
vanish up to order $2m$, and consequently we have that 
$\mathcal{T}_m^{(k)}(0)=0$ for $k=1,\dots,2m+1$. 
Since also $\mathcal{T}_m(0)=0$, \autoref{cor:symstep} implies that 
$\mathcal{T}_m\in\Sigma_{2m+1}$. 
\end{proof}

\begin{figure}[htbp]
    \centering
    \begin{subfigure}[t]{0.45\textwidth}
      \centering
      \includegraphics[width=0.9\textwidth]{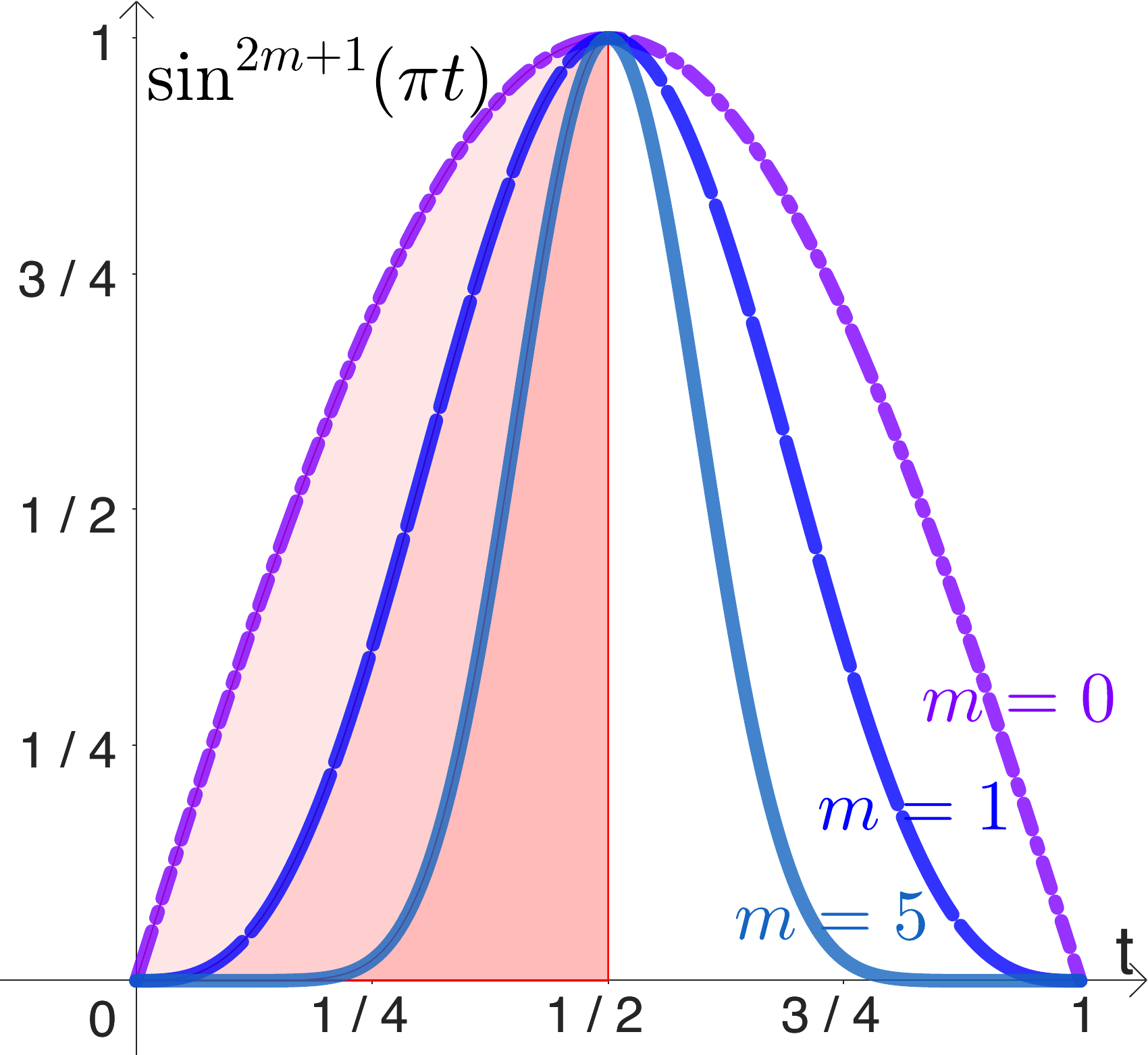}
      \caption{Odd powers of sine are bell-shaped with flat ends.}
    \end{subfigure}\hfill
    \begin{subfigure}[t]{0.45\textwidth}
      \centering
      \includegraphics[width=0.9\textwidth]{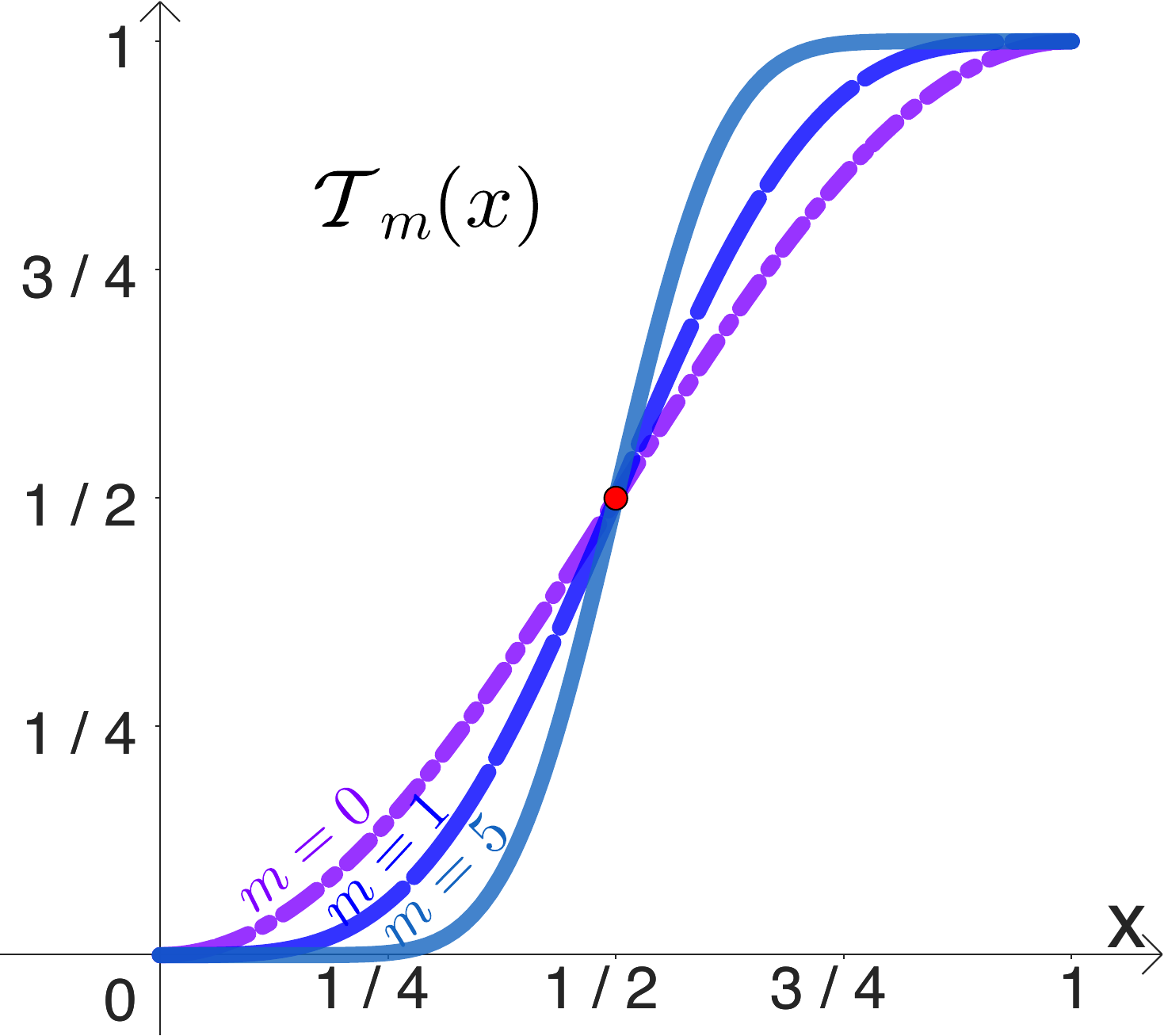}
      \caption{The trigonometric step function $\mathcal{T}_m$ has point symmetry.}
    \end{subfigure}
  \caption{The trigonometric step function $\mathcal{T}_{m}(x)$ is the
area under $\sin^{2m+1}(\pi t)$ from $t=0$ to $t=x$, scaled by the total area 
on $[0,1]$. As $m$ increases, $\sin^{2m+1}(\pi t)$ becomes flatter at the endpoints, 
while $\mathcal{T}_{m}$ becomes steeper near the midpoint.}
    \label{fig:trigstep}
\end{figure}

The smooth step function $\mathcal{T}_{m}$ can be thought of as a cumulative 
distribution function corresponding to a sine distribution supported on the unit
interval \cite{Edwards2000}, as depicted in \autoref{fig:trigstep}. Furthermore, 
it can be represented as a cosine expansion, described in what follows.

\begin{corollary}\label{cor:trigstep}
  Let $\mathcal{T}_m$ be as in \autoref{thm:trigstep}. Then
  \begin{equation}\label{eq:trigstep}
    \mathcal{T}_m(x) = \dfrac{1}{2}-\dfrac{1}{2a_m}
    \sum_{j=0}^m a_{m,j}\cos((2j+1)\pi x),\quad 0\leq x\leq 1,
  \end{equation}
  where 
  \begin{equation} a_{m,j} = (-1)^{m-j} \binom{2m+1}{m-j}/(2j+1), \qquad j=0,\dots,m
  \end{equation}
  and $a_m=a_{m,0}+\cdots+a_{m,m}$.
\end{corollary}
\begin{proof}
  Let $s_m(x):=\sum_{j=0}^m a_{m,j}\cos((2j+1)\pi x)$ for 
  $x\in[0,1]$. With the aid of \cite[p.~153]{Gradshteyn2007}, we have the antiderivative given by
  \begin{equation}
    \int\sin^{2m+1}(\pi x)dx = 
    \dfrac{(-1)^{m+1}}{2^{2m}\pi} \sum_{j=0}^m (-1)^{m-j} \binom{2m+1}{m-j}
    \dfrac{\cos((2j+1)\pi x)}{2j+1},
  \end{equation}
  and hence we deduce that 
  \begin{equation}
  S_m(x) = (-1)^{m+1}(s_m(x) - s_m(0))/(2^{2m}\pi),
  \end{equation}
  and so 
  \begin{equation}\mathcal{T}_m(x)=(s_m(x) - s_m(0))/(s_m(1) - s_m(0))
  \end{equation}
  for $0\leq x\leq 1$. Then \autoref{eq:trigstep} follows from 
  $s_m(0)=-s_m(1)=a_m$. 
\end{proof}

The coefficients of the trigonometric step function $\mathcal{T}_m$
satisfy a binomial identity.

\begin{corollary}\label{cor:binid}
  Let $m$ be a positive integer. Then
  \begin{equation}\label{eq:coeff}
  \sum_{j=0}^m (-1)^{m-j}(2j+1)^{2k-1}\binom{2m+1}{m-j} = 0,\quad k=1,\dots,m.
  \end{equation}
\end{corollary}
\begin{proof}
  Using induction and \autoref{cor:trigstep}, we can deduce that
  \begin{alignat}{3}
    \mathcal{T}^{(2k)}_m(x) && \ = \ & (-1)^{k+1}\dfrac{\pi^{2k}}{2a_m}
    \sum_{j=0}^m (2j+1)^{2k} a_{m,j} \cos((2j+1)\pi x), \quad
    && k=1,2,3,\dots. \label{eq:sigmaeven}
  \end{alignat}
  From \autoref{thm:trigstep}, we have $\mathcal{T}_m\in\Sigma_{2m+1}$, and
  consequently $\mathcal{T}_m^{(2k)}(0)=0$ for $k=1,\dots,m$.
  Then \autoref{eq:coeff} follows from evaluating \autoref{eq:sigmaeven} at $x=0$
  and substituting the values of $a_{m,j}$ given in \autoref{cor:trigstep}.
\end{proof}

The trigonometric step functions are related to a certain higher-order two-point boundary 
value problem. In order to reveal this relationship, we first consider the second-order
initial value problem given by the ordinary differential equation 
\begin{equation}\label{eq:ode2}
  (D^2+\pi^2I) y = \dfrac{\pi^2}{2},
\end{equation}
for $D^2 = d^2/dx^2$ and subject to the initial conditions
\begin{equation}\label{eq:ic2}
  y(0) = y'(0) = 0.
\end{equation}
The solution to \autoref{eq:ode2} and \autoref{eq:ic2} is given by 
\begin{equation}y_{0}(x)=\frac{1}{2}-\frac{1}{2}\cos(\pi x)\end{equation} for $x\geq 0$.
Then, from \autoref{cor:trigstep}, we deduce that $y_{0}(x)=\mathcal{T}_{0}(x)$ for $0\leq x\leq 1$, i.e., the trigonometric function $\mathcal{T}_{0}$ is a shifted harmonic oscillator 
restricted to $[0,1]$. Since $\mathcal{T}_{0}$ is a smooth step function of order 1 and agrees with 
$y_{0}$ on $[0,1]$, we deduce that $\mathcal{T}_{0}$ is the solution to the 
2nd-order two-point boundary value problem given by \autoref{eq:ode2} subject to 
\begin{equation}\label{eq:bc2}
  \begin{aligned}
    y(0) & = 0, & y(1) & = 1,\\
    y'(0) & = 0, & y'(1) & = 0.\\
  \end{aligned}
\end{equation}

Now, consider the order-$2m+2$ initial value problem given by
\begin{equation}\label{eq:odem}
  L_m[y] = (D^2+(2m+1)^2\pi^2I)\cdots(D^2+3^2\pi^2I)(D^2+\pi^2I) y = c_{m},
\end{equation}
subject to the initial conditions
\begin{equation}\label{eq:icm}
  y^{(k)}(0) = 0,\quad k=0,\dots,2m+1,
\end{equation}
where
\begin{equation}c_{m} = \dfrac{\left(1\cdot 3\cdots(2m+1)\pi^{m+1}\right)^{2}}{2}.\end{equation}

By the existence and uniqueness theorem for linear initial value problems,  
there is exactly one solution to \autoref{eq:odem} subject to \autoref{eq:icm}. 
In order to describe the homogeneous part of this solution, we substitute $y(x)=e^{\rho x}$ into the corresponding homogeneous equation, which gives the characteristic polynomial
\begin{equation}
  P(\rho)= (\rho^2+(2m+1)^2\pi^2)\cdots(\rho^2+3^2\pi^2)(\rho^2+\pi^2).
\end{equation}
The roots of $P$ are $\rho=\pm i(2j+1)\pi$ for $j=0,\dots,m$. Hence, the real-valued solution of the homogeneous equation corresponding to \autoref{eq:odem} can
be expressed as
\begin{equation}
  y_{H}(x) =
  \sum_{j=0}^{m} \alpha_{j}\cos((2j+1)\pi x)
  + \sum_{j=0}^{m} \beta_{j}\sin((2j+1)\pi x).
\end{equation}
Since a particular solution for \autoref{eq:odem} is $y_{P}(x)=1/2$, the general
solution is $y_m:=y_P+y_H$. The odd initial conditions in \autoref{eq:icm} give
\begin{equation}
\sum_{j=0}^{m}\beta_j(2j+1)^{2k+1}=0,\qquad k=0,\dots,m.
\end{equation}
This homogeneous linear system has only the trivial solution, and hence
$\beta_j=0$ for $j=0,\dots,m$. Therefore, any solution satisfying
\autoref{eq:odem} has the form given by
\begin{equation}
  y_m(x)=\frac{1}{2}+\sum_{j=0}^{m}\alpha_j\cos((2j+1)\pi x).
\end{equation}
From \autoref{cor:trigstep}, this is the same form as the trigonometric step function $\mathcal{T}_{m}$ with $\alpha_j = -a_{m,j}/(2a_m)$. By \autoref{thm:trigstep}, $\mathcal{T}_m$ satisfies the initial conditions given by \autoref{eq:odem}. Hence, from the existence and uniqueness theorem for linear initial value problems, it follows that $y = y_{m}(x)=\mathcal{T}_{m}(x)$ for $0\leq t\leq 1$ is the solution to the system given by \autoref{eq:odem} and \autoref{eq:icm}.


Given that $\mathcal{T}_{m}$ is a smooth step function of order $2m+1$ that agrees with 
$y_{m}$ on $[0,1]$, we have that $y_{m}(1)=1$ and $y_{m}^{(k)}(1) = 0$ for $k=1,\dots,2m+1$.
Hence, the trigonometric step function $\mathcal{T}_{m}$ is the solution to the order-$2m+2$ two-point boundary value 
problem given by \autoref{eq:odem} subject to
\begin{equation}\label{eq:bcm}
  \begin{aligned}
    y(0) & = 0, & y(1) & = 1, \\
    y'(0) & = 0, & y'(1) & = 0, \\
    & \vdots & & \vdots \\
    y^{(2m+1)}(0) & = 0, \quad & y^{(2m+1)}(1) & = 0. \\
  \end{aligned} 
\end{equation}

\begin{remark}
 Independently, the trigonometric smooth step function $\mathcal{T}_{m}$ was derived in \cite{Olofsen2019} in order to approximate the Fabius function, but without explicit expressions for the coefficients as in~\autoref{cor:trigstep}. There, it was shown that the coefficients $\alpha_j$ of $\mathcal{T}_m$ must sastify the equations given by
 \begin{alignat}{3}
 \sum_{j=0}^{m}\alpha_{j} & = -1/2, & \label{eq:coeff0} \\
 \sum_{j=0}^{m}(2j+1)^{2k}\alpha_{j} & = 0, & k=1,\dots,m \label{eq:coeffk}
 \end{alignat}
 in order for $T_m$ to be a smooth step function of order $2m+1$. We observe here that this is the same linear system one would obtain in order to satisfy the initial conditions of \autoref{eq:icm} upon substitution of $\mathcal{T}_m$ into ordinary differential equation of \autoref{eq:odem}. 

 To the best of our knowledge, neither the monotonicity property has been proved nor have explicit formulas for the coefficients been provided previously. As such, \autoref{thm:trigstep}, \autoref{cor:trigstep}, and \autoref{cor:binid} provide the corresponding results.


\end{remark}

\section*{Data availability statement}
This manuscript has no associated data.

\bibliography{bibliography} 

\section*{Declarations}

\subsection*{Funding}
This work was supported by the Agence Nationale de la Recherche (ANR) under Grant ANR-23-CE46-0008.

\subsection*{Competing interests}
The authors declare no conflicts of interest relevant to this work.

\end{document}